\documentclass{amsart}
\usepackage{amsthm,amsmath,amssymb,amscd,graphics,enumerate, stmaryrd,xspace,verbatim, epic, eepic,url}
\usepackage[all]{xypic}
\SelectTips{cm}{}

{
   \newtheorem{theorem}[subsubsection]{Theorem}
      \newtheorem*{theorem*}{Theorem}
   \newtheorem{proposition}[subsubsection]{Proposition}     
   \newtheorem{lemma}[subsubsection]{Lemma}

   \newtheorem{conjecture}[subsubsection]{Conjecture}
   \newtheorem*{conjecture*}{Conjecture}
   
}
{\theoremstyle{definition}
   \newtheorem{exercise}[subsubsection]{Exercise}       \newtheorem*{exercise*}{Exercise}
   \newtheorem{example}[subsubsection]{Example}
   \newtheorem*{example*}{Example}
   \newtheorem{definition}[subsubsection]{Definition}
   \newtheorem*{definition*}{Definition}
   \newtheorem{remark}[subsubsection]{Remark}
   
}
\newcommand{\RR}{{\mathbb{R}}}

\newcommand{\CC}{{\mathbb{C}}}
\newcommand{\QQ}{{\mathbb{Q}}}
\newcommand{\NN}{{\mathbb{N}}}
\newcommand{\MM}{{\mathbb{M}}}
\newcommand{\PP}{{\mathbb{P}}}
\newcommand{\ZZ}{{\mathbb{Z}}}

\renewcommand{\AA}{{\mathbb{A}}}

\newcommand{\bGamma}{{\boldsymbol{\Gamma}}}
\newcommand{\bDelta}{{\boldsymbol{\Delta}}}

\newcommand{\cC}{{\mathcal C}}

\newcommand{\cF}{{\mathcal F}}

\newcommand{\cK}{{\mathcal K}}

\newcommand{\cM}{{\mathcal M}}

\newcommand{\cO}{{\mathcal O}}
\newcommand{\cP}{{\mathcal P}}

\renewcommand{\cR}{{\mathcal R}}

\newcommand{\cX}{{\mathcal X}}
\newcommand{\cY}{{\mathcal Y}}

\def\<{\langle}
\def\>{\rangle}

\newcommand{\Spec}{\operatorname{Spec}}

\newcommand{\Pic}{{\operatorname{\mathbf{Pic}}}}

\newcommand{\Sym}{{\operatorname{Sym}}}

\newcommand{\das}{\dashrightarrow}

\newcommand{\ocC}{\overline{{\mathcal C}}}

\newcommand{\oC}{{\overline{C}}}

\newcommand{\oX}{{\overline{X}}}

\newcommand{\oSigma}{{\overline{\Sigma}}}
\newcommand{\OKS}{{\cO_{k,S}}}
\newcommand{\double}{\genfrac..{0pt}1
{\raise -1pt\hbox{$\scriptstyle\longrightarrow$}}{\raise 3pt\hbox
{$\scriptstyle\longrightarrow$}}} 

\newcommand{\setmin}{\smallsetminus}

\renewcommand{\setminus}{\smallsetminus}

\def\tototi{\mathbin{\mathop{\otimes}\limits^{\raise-1pt\hbox
{$\scriptscriptstyle {\rm L}$}}}}

\def\indlim{\mathop{\vrule width0pt height7pt depth
4pt\smash{\lim\limits_{\raise 1pt\hbox to 14.5pt
{\rightarrowfill}}}}}
\def\projlim{\mathop{\vrule width0pt height7pt depth
4pt\smash{\lim\limits_{\raise 1pt\hbox to 14.5pt
{\leftarrowfill}}}}}

\newcommand\displaceamount{3pt}

\newcommand{\doubledown}{\ar@<\displaceamount>[d]\ar@<-\displaceamount>[d]}

\newcommand{\doubleup}{\ar@<\displaceamount>[u]\ar@<-\displaceamount>[u]}

\newcommand{\doubleright}{\ar@<\displaceamount>[r]\ar@<-\displaceamount>[r]}

\newcommand{\lecture}[1]{\section{#1}}

\begin{document}

\title{Birational geometry for number theorists}

\author[D. Abramovich]{Dan Abramovich}
\address{Department of Mathematics, Box 1917, Brown University,
Providence, RI, 02912, U.S.A} 
\email{abrmovic@math.brown.edu}
\thanks{Partial support for research provided by NSF grants DMS-0301695 and DMS-0603284.}


\maketitle

\setcounter{tocdepth}{1}
\tableofcontents

\section*{Introduction} 
\renewcommand{\thesection}{Lecture \arabic{section}}
\renewcommand{\thesubsection}{\arabic{section}.\arabic{subsection}}

When thinking about the course ``birational geometry for number theorists" I so na\"{\i}vely agreed to give at the G\"ottingen summer school, I could not avoid imagining the spirit of the late Serge Lang, not so quietly beseeching one to do things right, keeping the theorems functorial with respect to ideas, and definitions natural. But most important is the fundamental tenet of diophantine geometry,  for which Lang was one of the strongest and loudest advocates, which was so aptly summarized in the introduction of Hindry-Silverman \cite{Hindry-Silverman}:
\begin{center}
\fbox{\large \bf GEOMETRY DETERMINES ARITHMETIC.}
\end{center}

To make sense of  this, largely conjectural,  epithet, it is good to have some loose background in birational geometry, which I will try to provide. For the arithmetic motivation I will explain conjectures of Bombieri, Lang and Vojta, and new and exciting versions of those due to Campana. In fact, I imagine Lang would insist (strongly, as only he could) that Campana's conjectures most urgently need further investigation, and indeed in some sense they form the centerpiece of these notes.

Birational geometry is undergoing revolutionary developments these very days: large portions of the minimal model program were solved soon after the G\"ottingen lectures \cite{BCHM}, and it seems likely that more is to come. Also a number of people seem to have made new inroads into the long standing resolution of singularities problem.  I am not able to report on the latter, but I will give a brief account of the minimal model program as it seems to stand at this point in time.

Our convention: a variety over $k$ is an {\em absolutely} reduced and irreducible  scheme of finite type over $k$.

{\sc Acknowledgements:} I thank the CMI and the organizers for
inviting me, I thank the colleagues and students at Brown for their
patience with my ill prepared preliminary lectures and numerous
suggestions,  I thank F. Campana for a number of inspiring
discussions, H.-H. Tseng for a number of good comments, and L. Caporaso for the notes of her MSRI lecture
\cite{Caporaso}, to which my lecture plans grew increasingly
close. The treatment of the minimal model program is influenced by
lectures of Ch. Hacon and J. McKernan and discussions with them. Many thanks are due to the referee who caught a large number of errors and made numerous suggestions. Of course all remaining errors are entirely my responsibility. Anything new is partially supported by the NSF grants DMS-0301695 and DMS-0603284. 

\setcounter{section}{-1}

\lecture{Geometry and arithmetic of curves} 

The arithmetic of algebraic curves is one area where basic relationships between geometry and arithmetic are known, rather than conjectured. Much of the material here is covered in Darmon's lectures of this summer school.

\subsection{Closed curves}

Consider a smooth projective algebraic curve $C$ defined over a number field $k$. We are interested in a qualitative relationship between its arithmetic and geometric properties.

We have three basic facts:

\subsubsection{} A curve of genus 0 becomes rational after at most a quadratic extension $k'$ of $k$, in which case its set of rational points $C(k')$ is infinite (and therefore dense in the Zariski topology). 

\subsubsection{} A curve of genus 1 has a rational point after a finite extension $k'$ of $k$ (though the degree is not a-priori bounded), and has positive Mordell--Weil rank after a further quadratic extension $k''/k'$, in which case again its set of rational points $C(k'')$ is infinite (and therefore dense in the Zariski topology). 

We can immediately introduce the following definition:
\begin{definition}\label{Def:pot-dense}
Let $X$ be an algebraic variety defined over $k$. We say that rational points on $X$ are potentially dense, if there is a finite extension $k'/k$ such that the set $X(k')$ is dense in $X_{k'}$ in the Zariski topology.
\end{definition}

Thus rational points on a curve of genus 0 or 1 are potentially dense.

Finally we have
\begin{theorem}[Faltings, 1983]
Let $C$ be an algebraic curve of genus $>1$ over a number field $k$. Then $C(k)$ is finite.
\end{theorem}

See, e.g. \cite{Faltings, Hindry-Silverman}.

In other words, rational points on a curve $C$ of genus $g$ are potentially dense if and only if $g\leq 1$.

\subsubsection{} So far there isn't much birational geometry involved, because we have the old theorem:
\begin{theorem}
A smooth algebraic curve is uniquely determined by its function field.
\end{theorem} 

But this is an opportunity to introduce a tool: on the curve $C$ we have a canonical divisor class $K_C$, such that $\cO_C(K_C) = \Omega^1_C$, the sheaf of differentials, also known by the notation $\omega_C$ - the dualizing sheaf. We have:
\begin{enumerate}
\item $\deg K_C = 2g-2 = -\chi(C_\CC)$, where $\chi(C_\CC)$ is the topological Euler characteristic of the complex Riemann surface $C_\CC$.
\item $\dim H^0(C, \cO_C(K_C)) = g$.
\end{enumerate}

For future discussion, the first property will be useful. We can now summarize, following \cite{Hindry-Silverman}: 

\subsubsection{} \label{Tab:rational}
\begin{center}\large  \begin{tabular}{|l|l|}\hline 
 Degree of $K_C$ & rational points\\ \hline\hline
$2g-2\leq 0$ & potentially dense \\ \hline
$2g-2> 0$ & finite \\ \hline
\end{tabular} \end{center}

\subsection{Open curves} 

\subsubsection{}

Consider a smooth quasi-projective algebraic curve $C$ defined over a number field $k$. It has a unique smooth projective completion $C \subset \oC$, and the complement is a finite set $\Sigma = \oC \setminus C$. Thinking of $\Sigma$ as a reduced divisor of some degree $n$, a natural line bundle to consider is $\cO_{\oC}(K_\oC + \Sigma)\simeq \omega_{\cO}(\Sigma)$, the sheaf of differentials with logarithmic poles on $\Sigma$, whose degree is again $-\chi^{top}(C) = 2g-2+n$. The sign of $2g-2+n$ again serves as the geometric invariant to consider.

\subsubsection{} Consider for example   the affine line. Rational points on the affine line are not much more interesting than those on $\PP^1$. But we can also consider the behavior of {\em integral} points, where interesting results do arise. However, what does one mean by integral points on $\AA^1$?
The key is that integral points are an invariant of an ``integral  model" of $\AA^1$ over $\ZZ$.



\subsubsection{}
Consider the ring of integers $\cO_k$ and a finite set $S \subset \Spec \cO_k$ of finite primes. One can associate to it the ring $\OKS$ of $S$-integers, of elements in $K$ which are in $\cO_\wp$ for any prime $\wp \not\in S$. 

Now consider a {\em model} of $C$ over $\cO_{k,S}$, namely a scheme $\cC$ of finite type over $\cO_{k,S}$ with an isomorphism of the generic fiber $\cC_k \simeq C$. It is often useful to start with a model $\ocC$ of $\oC$, and take $\cC = \ocC \setminus \oSigma$.

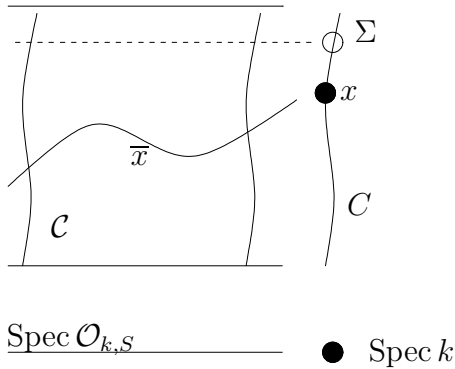
\begin{center}
\setlength{\unitlength}{0.0005in}
\begingroup\makeatletter\ifx\SetFigFont\undefined%
\gdef\SetFigFont#1#2#3#4#5{%
  \reset@font\fontsize{#1}{#2pt}%
  \fontfamily{#3}\fontseries{#4}\fontshape{#5}%
  \selectfont}%
\fi\endgroup%
{\renewcommand{\dashlinestretch}{30}
\begin{picture}(4285,3906)(0,-10)
\put(3312,2979){\blacken\ellipse{212}{212}}
\put(3312,2979){\ellipse{212}{212}}
\put(3387,279){\blacken\ellipse{212}{212}}
\put(3387,279){\ellipse{212}{212}}
\path(12,3879)(2862,3879)
\drawline(2862,3879)(2862,3879)
\path(12,1179)(2862,1179)
\path(12,279)(2862,279)
\path(3462,3804)(3462,3802)(3460,3797)
	(3458,3787)(3455,3773)(3451,3753)
	(3446,3727)(3439,3695)(3431,3658)
	(3423,3617)(3414,3573)(3405,3526)
	(3395,3478)(3386,3429)(3377,3381)
	(3368,3334)(3360,3288)(3353,3244)
	(3346,3202)(3339,3162)(3334,3124)
	(3329,3088)(3325,3054)(3321,3021)
	(3318,2989)(3316,2958)(3314,2928)
	(3313,2899)(3312,2870)(3312,2841)
	(3312,2811)(3313,2780)(3314,2749)
	(3316,2718)(3318,2686)(3320,2654)
	(3323,2621)(3327,2588)(3330,2555)
	(3334,2521)(3338,2487)(3342,2453)
	(3347,2418)(3351,2384)(3356,2350)
	(3360,2316)(3365,2283)(3369,2250)
	(3374,2217)(3378,2185)(3382,2154)
	(3385,2124)(3388,2094)(3391,2065)
	(3394,2036)(3396,2009)(3398,1981)
	(3400,1954)(3401,1925)(3401,1896)
	(3401,1867)(3400,1837)(3399,1807)
	(3397,1776)(3395,1743)(3392,1709)
	(3388,1673)(3384,1635)(3378,1595)
	(3373,1553)(3367,1510)(3360,1466)
	(3353,1422)(3346,1379)(3339,1338)
	(3333,1300)(3327,1266)(3322,1237)
	(3318,1215)(3315,1198)(3314,1188)
	(3312,1182)(3312,1179)
\path(2637,3804)(2637,3802)(2635,3797)
	(2633,3787)(2630,3773)(2626,3753)
	(2621,3727)(2614,3695)(2606,3658)
	(2598,3617)(2589,3573)(2580,3526)
	(2570,3478)(2561,3429)(2552,3381)
	(2543,3334)(2535,3288)(2528,3244)
	(2521,3202)(2514,3162)(2509,3124)
	(2504,3088)(2500,3054)(2496,3021)
	(2493,2989)(2491,2958)(2489,2928)
	(2488,2899)(2487,2870)(2487,2841)
	(2487,2811)(2488,2780)(2489,2749)
	(2491,2718)(2493,2686)(2495,2654)
	(2498,2621)(2502,2588)(2505,2555)
	(2509,2521)(2513,2487)(2517,2453)
	(2522,2418)(2526,2384)(2531,2350)
	(2535,2316)(2540,2283)(2544,2250)
	(2549,2217)(2553,2185)(2557,2154)
	(2560,2124)(2563,2094)(2566,2065)
	(2569,2036)(2571,2009)(2573,1981)
	(2575,1954)(2576,1925)(2576,1896)
	(2576,1867)(2575,1837)(2574,1807)
	(2572,1776)(2570,1743)(2567,1709)
	(2563,1673)(2559,1635)(2553,1595)
	(2548,1553)(2542,1510)(2535,1466)
	(2528,1422)(2521,1379)(2514,1338)
	(2508,1300)(2502,1266)(2497,1237)
	(2493,1215)(2490,1198)(2489,1188)
	(2487,1182)(2487,1179)
\path(312,3804)(312,3802)(310,3797)
	(308,3787)(305,3773)(301,3753)
	(296,3727)(289,3695)(281,3658)
	(273,3617)(264,3573)(255,3526)
	(245,3478)(236,3429)(227,3381)
	(218,3334)(210,3288)(203,3244)
	(196,3202)(189,3162)(184,3124)
	(179,3088)(175,3054)(171,3021)
	(168,2989)(166,2958)(164,2928)
	(163,2899)(162,2870)(162,2841)
	(162,2811)(163,2780)(164,2749)
	(166,2718)(168,2686)(170,2654)
	(173,2621)(177,2588)(180,2555)
	(184,2521)(188,2487)(192,2453)
	(197,2418)(201,2384)(206,2350)
	(210,2316)(215,2283)(219,2250)
	(224,2217)(228,2185)(232,2154)
	(235,2124)(238,2094)(241,2065)
	(244,2036)(246,2009)(248,1981)
	(250,1954)(251,1925)(251,1896)
	(251,1867)(250,1837)(249,1807)
	(247,1776)(245,1743)(242,1709)
	(238,1673)(234,1635)(228,1595)
	(223,1553)(217,1510)(210,1466)
	(203,1422)(196,1379)(189,1338)
	(183,1300)(177,1266)(172,1237)
	(168,1215)(165,1198)(164,1188)
	(162,1182)(162,1179)
\path(12,2004)(14,2006)(18,2010)
	(26,2017)(38,2028)(55,2044)
	(77,2064)(104,2089)(136,2117)
	(172,2150)(211,2185)(254,2222)
	(298,2260)(343,2298)(388,2336)
	(433,2373)(477,2408)(519,2441)
	(560,2473)(599,2501)(637,2527)
	(673,2551)(707,2571)(739,2590)
	(770,2606)(800,2619)(829,2630)
	(857,2639)(884,2646)(910,2650)
	(936,2653)(962,2654)(989,2653)
	(1017,2650)(1045,2646)(1072,2639)
	(1101,2632)(1129,2622)(1158,2612)
	(1187,2599)(1217,2586)(1248,2572)
	(1278,2557)(1310,2541)(1341,2524)
	(1373,2507)(1405,2490)(1437,2472)
	(1469,2455)(1501,2439)(1533,2422)
	(1565,2407)(1597,2392)(1628,2378)
	(1660,2365)(1691,2354)(1721,2344)
	(1752,2335)(1783,2328)(1813,2322)
	(1844,2319)(1875,2316)(1902,2316)
	(1930,2318)(1958,2320)(1987,2325)
	(2017,2332)(2048,2340)(2081,2350)
	(2114,2362)(2150,2377)(2187,2393)
	(2226,2411)(2266,2432)(2309,2455)
	(2354,2480)(2401,2507)(2450,2535)
	(2500,2566)(2551,2597)(2603,2630)
	(2654,2663)(2705,2696)(2754,2728)
	(2800,2759)(2843,2788)(2882,2814)
	(2916,2837)(2944,2857)(2968,2873)
	(2985,2885)(2998,2894)(3006,2900)
	(3010,2903)(3012,2904)
\dashline{60.000}(3162,3504)(87,3504)
\put(3387,3504){\ellipse{212}{212}}
\put(3762,129){\makebox(0,0)[lb]{{\SetFigFont{12}{14.4}{\rmdefault}{\mddefault}{\updefault}$\Spec k$}}}
\put(0,254){\makebox(0,0)[lb]{{\SetFigFont{12}{14.4}{\rmdefault}{\mddefault}{\updefault}$\Spec \OKS$}}}
\put(3612,3504){\makebox(0,0)[lb]{{\SetFigFont{12}{14.4}{\rmdefault}{\mddefault}{\updefault}$\Sigma$}}}
\put(3462,2904){\makebox(0,0)[lb]{{\SetFigFont{12}{14.4}{\rmdefault}{\mddefault}{\updefault}$x$}}}
\put(1287,2229){\makebox(0,0)[lb]{{\SetFigFont{12}{14.4}{\rmdefault}{\mddefault}{\updefault}$\overline
      x$}}}
\put(3537,1704){\makebox(0,0)[lb]{{\SetFigFont{12}{14.4}{\rmdefault}{\mddefault}{\updefault}$C$}}}
\put(462,1479){\makebox(0,0)[lb]{{\SetFigFont{12}{14.4}{\rmdefault}{\mddefault}{\updefault}$\cC$}}}
\end{picture}
}
\end{center}

Now it is clear how to define integral points:  an $S$-integral point on $\cC$ is  simply an element of $\cC(\cO_{k,S})$, in other words, a section of $\cC \to \Spec (\cO_{k,S})$.  This is related to rational points on a proper curve as follows:

\subsubsection{} If $\Sigma = \emptyset$, and the model chosen is proper, the notions of integral and rational points agree, because of the valuative criterion for properness.

\begin{exercise} Prove this!
\end{exercise}

 We have the following facts:
\subsubsection{} If $C$ is rational and $n\leq 2$, then after possibly enlarging $k$ and $S$, any integral model of $C$ has an infinite collection of integral points.

\begin{exercise} Prove this!
\end{exercise}

On the other hand, we have:
\begin{theorem}[Siegel's Theorem] If  $n\geq 3$, or if $g>0$ and $n>0$, then for any integral model $\cC$ of $C$, the set of integral points  $\cC(\OKS)$ is finite.
\end{theorem}

A good generalization of Definition \ref{Def:pot-dense} is the following:
\begin{definition}
Let $X$ be an algebraic variety defined over $k$ with a model $\cX$ over $\OKS$. We say that integral points on $X$ are potentially dense, if there is a finite extension $k'/k$, and an enlargement $S'$ of the set of places in $k'$ over $S$, such that the set $\cX(\cO_{k',S'})$ is dense in $X_{k'}$ in the Zariski topology.
\end{definition}

We can apply this definition in the case of a curve $C$ and generalize \ref{Tab:rational}, as in \cite{Hindry-Silverman}, as follows:

\subsubsection{} \label{Tab:integral}
\begin{center} \large \begin{tabular}{|l|l|}\hline 
degree of $K_\oC + \Sigma$ & integral points\\ \hline\hline
$2g-2+n\leq 0$ & potentially dense \\ \hline
$2g-2+n> 0$ & finite \\ \hline
\end{tabular}
\end{center}

\subsubsection{}\label{Lesson:open} One lesson we must remember from this discussion is that  \begin{center}
\fbox{\large For {\bf open} varieties we use {\bf integral} points on {\bf integral models}.}
\end{center}

\subsection{Faltings implies Siegel}

Siegel's theorem was proven years before Faltings's theorem. Yet it is instructive, especially in the later parts of these notes, to give the following argument showing that Faltings's theorem implies Siegel's. 

\begin{theorem}[Hermite-Minkowski, see \cite{Hindry-Silverman} page 264]
Let $k$ be a number field, $S\subset \Spec \OKS$ a finite set of finite places, and $d$ a positive integer. Then there are only finitely many extensions $k'/k$ of degree $\leq d$ unramified outside $S$.
\end{theorem}

From which one can deduce

\begin{theorem}[Chevalley-Weil, see \cite{Hindry-Silverman} page 292]
Let $\pi:\cX \to \cY$ be a finite \'etale morphism of schemes over $\OKS$. Then there is a finite extension $k'/k$, with $S'$ lying over $S$, such that $\pi^{-1}\cY(\OKS) \subset \cX(\cO_{k',S'})$. 
\end{theorem}

On the geometric side we have an old topological result

\begin{theorem}
If $C$ is an open curve with $2g-2+n>0$ and $n>0$, defined over $k$, there is a finite extension $k'/k$ and a finite unramified covering $D \to C$, such that $g(D)>1$.
\end{theorem}

\begin{exercise}
Combine these theorems to obtain a proof of Siegel's theorem assuming Faltings's theorem.
\end{exercise}

This is discussed in Darmon's lectures, as well as \cite{Hindry-Silverman}.

\subsubsection{} Our lesson this time  is that
  \begin{center}
\fbox{\large Rational and integral points can be controlled in finite \'etale covers.}
\end{center}

\subsection{Function field case}
There is an old and distinguished tradition of comparing results over number fields with results over function fields. To avoid complications I will concentrate on function fields of characteristic 0, and consider closed curves only. 

\subsubsection{} If $K$ is the function field of a complex variety $B$, then a variety $X/K$ is the generic fiber of a scheme $\cX/B$, and  a $K$-rational point $P\in X(K)$ can be thought of as a rational section of $\cX \to B$.   If $B$ is a smooth curve and $\cX \to B$ is proper, then again a $K$-rational point $P\in X(K)$ is equivalent to a {\em regular} section $B \to \cX$.

\begin{exercise} Make sense of this (i.e. prove this)! 
\end{exercise}

\subsubsection{} The notion of {\em integral points} is similarly defined  using sections. When $\dim B>1$ there is an intermediate notion of {\em properly rational points}: a $K$-rational point $p$ of $X$ is a properly rational point of $\cX/B$  if the closure $B'$ of $p$ in $\cX$ maps properly to $B$.

Consider now $C/K$ a curve. Of course it is possible that  $\cC$ is, or is birationally equivalent to, $C_0\times B$, in which case we have plenty of constant sections coming from $C_0(\CC)$, corresponding to constant points in $C(K)$. But that is almost all there is:

\begin{theorem}[Manin \cite{Manin}, Grauert \cite{Grauert}] let $k$ be a field of characteristic 0, let $K$ be a regular extension of $k$, and $C/K$ a smooth curve. 
Assume $g(C)>1$. If $C(K)$ is infinite, then there is a curve $C_0/k$ with $(C_0)_K \simeq C$, such that $C(K)\setminus C_0(k)$ is finite.
\end{theorem}

\begin{exercise}
What does this mean for constant curves $C_0\times B$ in terms of maps from $C_0$ to $B$?
\end{exercise}

Working inductively on transcendence degree, and using Faltings's Theorem, we obtain:
\begin{theorem} Let $C$ be a curve of genus $>1$ over a field $k$ finitely generated over $\QQ$. 
Then the set of $k$-rational points $C(k)$ is finite.
\end{theorem}
\begin{exercise}
Prove this, using previous results as given!
\end{exercise}
See \cite{Samuel}, \cite{Moosa-Scanlon} for  appropriate statements in positive characteristics.

\lecture{Kodaira dimension}

\subsection{Iitaka dimension}
Consider now a smooth, projective variety $X$ of dimension $d$ over a field $k$ of characteristic 0.
We seek an analogue of the sign of $2g-2$ in this case. The approach is by counting sections of the canonical line bundle $\cO_X(K_X) = \wedge^d\Omega^1_X$. Iitaka's book \cite{Iitaka} is a good reference.

\begin{theorem}\label{Th:Iitaka-dim}
Let  $L$ be a line bundle on $X$. Assume $h^0(X, L^n)$ does not vanish for all positive integers $n$. Then there is a unique integer $\kappa = \kappa(X, L)$ with $0\leq \kappa\leq d$ such that 
$$\limsup\limits_{n\to \infty} \frac{h^0(X, L^n)}{n^\kappa}$$ exists and is nonzero. 
\end{theorem}

\begin{definition}
\begin{enumerate}
\item 
The integer $\kappa(X, L)$ in the theorem is called {\em the Iitaka dimension of $(X,L)$}.
\item  In the special case $L = \cO_X(K_K)$ we write $\kappa(X):= \kappa(X, L)$ and call $\kappa(X)$ {\em the Kodaira dimension of $X$}.
\item 
 It is customary to set $\kappa(X,L)$ to be either $-1$ or  $-\infty$ if $h^0(X, L^n)$ vanishes  for all positive integers $n$. It is safest to say in this case that the Iitaka dimension is {\em negative}. I will use $-\infty$.
 \end{enumerate}
 \end{definition}
 
 We will see an algebraic justification for the  $-1$ convention immediately in  Proposition \ref{prop:kod-tr-deg},  and a geometric justification for the more commonly used $-\infty$ in a bit.

An algebraically meaningful presentation of the Iitaka dimension is the following: 

\begin{proposition}\label{prop:kod-tr-deg} Consider the algebra of sections 
$$\cR(X,L) := \bigoplus\limits_{n\geq 0} H^0(X, L^n).$$
Then, with the $-1$ convention, 
$$\operatorname{tr.deg} \cR(X,L) = \kappa(X,L) + 1.$$
\end{proposition}

\begin{definition}
We say that a property holds for a sufficiently high and divisible $n$, if there exists $n_0>0$ such that the property holds for every positive multiple of $n_0$.
\end{definition}

A geometric meaning of $\kappa(X,L)$ is given by the following: 
\begin{proposition} Assume $\kappa(X,L)\geq 0$.  Then for sufficiently high and divisible $n$, the dimension of the image  of the rational map $\phi_{L^n}: X \das \PP H^0(X, L^n)$ is precisely $\kappa(X,L)$.
\end{proposition}

Even more precise is:

\begin{proposition} \label{Prop:Iitaka-image} There is $n_0>0$ such that the image $\phi_{L^n}(X)$  is birational to $\phi_{L^{n_0}}(X)$ for all $n>0$ divisible by $n_0$. 
\end{proposition}

\begin{definition}
\begin{enumerate}
\item  The birational equivalence class of the variety $\phi_{L^{n_0}}(X)$ is denoted $I(X,L)$. 
\item The rational map $X \to I(X,L)$ is called the {\em Iitaka fibration} of $(X,L)$. 
\item In case $L$ is the canonical bundle $\omega_X$,  this map is simply called  the Iitaka fibration of $X$, written $X \to I(X)$ 
\end{enumerate}
\end{definition}
 
The following notion is important:
\begin{definition}
The variety $X$ is said to be {\em of general type} of $\kappa(X) = \dim X$.
\end{definition}
\begin{remark} The name is not as informative as one could wish. It comes from the observation that surfaces not of general type can be nicely classified, whereas there is a whole zoo of surfaces of general type.
\end{remark}

\begin{exercise}
Prove Proposition \ref{Prop:Iitaka-image}:
\begin{enumerate}
\item Show that if $n,d>0$ and $H^0(X, L^n) \neq 0$ then there is a dominant rational map $\phi_{L^{nd}}(X ) \das \phi_{L^n}(X)$ such that the following diagram is commutative:
$$\xymatrix{
X \ar@{-->}^{\phi_{L^{nd}}}[r]\ar@{-->}_{ \phi_{L^n}}[rd] & \phi_{L^{nd}}(X )\ar@{-->}[d] \\ & \phi_{L^n}(X).
}$$
\item Conclude that  $\dim  \phi_{L^n}(X)$ is a constant $\kappa$ for large and divisible $n$. 
\item Suppose $n>0$ satisfies $\kappa:=\dim  \phi_{L^n}(X)$. Show that for any $d>0$, the function field of $\phi_{L^{nd}}(X )$ is algebraic over the function field $\phi_{L^n}(X)$.
\item Recall that for any variety $X$, any subfield $L$ of $K(X)$ containing $k$ is finitely generated. Apply this to the algebraic closure of $ \phi_{L^n}(X)$ in $K(X)$ to complete the proof of the proposition.
\end{enumerate}
\end{exercise}

For details see \cite{Iitaka}.

\begin{exercise}
Use proposition \ref{Prop:Iitaka-image} to prove Theorem \ref{Th:Iitaka-dim}.
\end{exercise}

\subsection{Properties and examples of the Kodaira dimension}

\begin{exercise}
Show that $\kappa(\PP^n) = -\infty$ and $\kappa(A) = 0$ for an abelian variety $A$.
\end{exercise}

\subsubsection{Curves:}

\begin{exercise*}
Let $C$ be a smooth projective curve and $L$ a line bundle.  Prove that 
$$\kappa(C, L) = \begin{cases}
1 & \text{if  $\deg_CL>0,$} \\
0 & \text{if  $L$  is torsion, and } \\
<0& \text{otherwise.}
\end{cases}$$
In particular, 
$$\kappa(C) = \begin{cases}
1 & \text{if  $g>1$,} \\
0 & \text{if  $g=1$, and } \\
<0& \text{if  $g=0$.}
\end{cases}$$
\end{exercise*}

\subsubsection{Birational invariance:}
\begin{exercise*}
Let $X' \das X$ be a birational map of smooth projective varieties. 
Show that the spaces $H^0(X, \cO_X(mK_X))$ and $H^0(X', \cO_{X'}(mK_{X'}))$ are canonically isomorphic.
Deduce that $\kappa(X) = \kappa(X')$.
\end{exercise*}
 (See \cite{Hartshorne}, Chapter II, Theorem 8.19).

\subsubsection{Generically finite dominant maps}
\begin{exercise*}
Let $f:X' \to X$ be a generically finite map of smooth projective varieties.

Show that $\kappa(X') \geq \kappa(X)$. 
\end{exercise*}
\subsubsection{Finite \'etale maps}
\begin{exercise*}
Let $f:X' \to X$ be a finite \'etale map of smooth projective varieties. 

Show that $\kappa(X') = \kappa(X)$.
\end{exercise*}

\subsubsection{Field extensions:}
\begin{exercise*}
Let $k'/k$ be a field extension, $X$ a variety over $k$ with line bundle $L$, and $X_{k'}, L_{k'}$ the result of base change. 

Show that  $\kappa(X,L) = \kappa(X_{k'}, L_{k'})$. In particular $\kappa(X) = \kappa(X_{k'})$.
\end{exercise*}

\subsubsection{Products}
\begin{exercise*} Show that, with the $-\infty$ convention,
$$\kappa(X_1\times X_2, L_1\boxtimes L_2) = \kappa(X_1,L_1) + \kappa(X_2,L_2).$$ 

Deduce that $\kappa(X_1\times X_2) = \kappa(X_1) + \kappa(X_2)$. 
\end{exercise*}

This so called ``easy additivity" of the Kodaira dimension  is the main reason for the $-\infty$ convention.

\subsubsection{Fibrations}
The following is subtle and difficult:

\begin{theorem*}[Siu's theorem on deformation invariance of plurigenera \cite{Siu,Siu1}]
Let $X \to B$ be a smooth projective morphism with connected geometric fibers, and $m$ a positive integer. Then for closed points $b\in B$, the dimension $h^0(X_b, \cO(mK_{X_b}))$ is independent of $b\in B$.  In particular $\kappa(X_b)$ is independent of the closed point $b\in B$.
\end{theorem*}

\begin{exercise} Let $X \to B$ be a morphism of smooth projective varieties with connected geometric fibers. Let $b\in B$ be such that $X\to B$ is smooth over $b$, and let $\eta_B\in B$ be the generic point.

Use ``cohomology and base change" and Siu's theorem  to deduce that 
$$\kappa(X_b) = \kappa(X_{\eta_B}).$$ 
\end{exercise}

\begin{definition} The Kodaira defect of $X$ is $\delta(X) = \dim(X) - \kappa(X)$. 
\end{definition}

\begin{exercise} 
 Let $X \to B$ be a morphism of smooth projective varieties with connected geometric fibers. Show that the Kodaira defects satisfy $\delta(X) \geq \delta(X_{\eta_B})$. Equivalently $\kappa(X) \leq  \dim(B) + \kappa(X_{\eta_B})$.
 \end{exercise}
 
We remark that  before Siu's deformation invariance theorem was proven, a weaker and more technical result, yet still very useful, saying that the Kodaira dimension is constant  on ``very general fibers" was used.

 \begin{exercise} \label{Ex:easy-inequality}
 Let $Y \to B$ be a morphism of smooth projective varieties with connected  geometric fibers, and $Y\to X$ a generically finite map. Show that $\delta(X) \geq \delta(Y_{\eta_B})$. In other words, $\kappa(X) \leq \kappa(Y_{\eta_B}) + \dim B$. 
 \end{exercise}
 
 This so called ``easy subadditivity" has many useful consequences.
 
 \begin{definition} We say that $X$ is uniruled if 
there is a variety $B$ of dimension $\dim X - 1$ and a dominant rational map $B\times \PP^1 \das X$.  
 \end{definition}
 
 \begin{exercise}
 If $X$ is uniruled, show that $\kappa(X) = -\infty$. 
 \end{exercise}

The converse is an important conjecture, sometimes known as the $(-\infty)$-Conjecture. It is a consequence of the ``good minimal model" conjecture:

\begin{conjecture}\label{Conj:k<0-uniruled}
Assume $X$ is not uniruled. Then $\kappa(X) \geq 0$.
\end{conjecture}

 \begin{exercise}
 If $X$ is covered by a family of elliptic curves, show that $\kappa(X) \leq \dim X - 1$. 
 \end{exercise}

\subsubsection{Surfaces}

Surfaces of Kodaira dimension $<2$ are ``completely classified". Some of these you can place in the table using what you have learned so far. In the following description we give a representative of the birational class of each type:

\begin{center}\begin{tabular}{|l|l|}\hline 
$\kappa$ & description \\ \hline \hline 
$-\infty$ & ruled surfaces: $\PP^2$ or $\PP^1 \times C$ \\ \hline 
$0$ & a. abelian surfaces \\
 & b. bielliptic surfaces\\ 
 & k. K3 surfaces\\
 & e. Enriques surfaces\\
 \hline
 $1$ & all other elliptic surfaces\\ \hline
\end{tabular}\end{center}

\subsubsection{Iitaka's program}

Here is a central conjecture of birational geometry:
\begin{conjecture*}[Iitaka]
Let $X\to B$ be a surjective morphism of smooth projective varieties. Then 
$$\kappa(X) \geq \kappa(B) + \kappa(X_{\eta_B}).$$
\end{conjecture*}
\subsubsection{}
Major progress on this conjecture was made through the years by several geometers, including  Fujita \cite{Fujita}, Kawamata \cite{Kawamata}, Viehweg \cite{Viehweg-additivity} and Koll\'ar \cite{Kollar-subadditivity}. The key, which makes this conjecture plausible, is the semi-positivity properties of   the relative dualizing sheaf $\omega_{X/B}$, which originate from work of Arakelov and rely on deep Hodge theoretic arguments. 

Two results will be important for these lectures.
\begin{theorem}[Kawamata]
Iitaka's conjecture follows from the Minimal Model Program: if $X_{\eta_B}$ has a good minimal model then $\kappa(X) \geq \kappa(B) + \kappa(X_{\eta_B}).$
\end{theorem}

\begin{theorem}[Viehweg]
Iitaka's conjecture holds in case $B$ is of general type, namely:

Let $X\to B$ be a surjective morphism of smooth projective varieties, and assume $\kappa(B) = \dim B$. Then 
$\kappa(X) = \dim(B) + \kappa(X_{\eta_B}).$
\end{theorem}

Note that equality here is forced by the easy subadditivity inequality: $\kappa(X) \leq  \dim(B) + \kappa(X_{\eta_B})$ always holds.

\begin{exercise} Let $X, B_1, B_2$ be smooth projective varieties.
Suppose $X \to B_1\times B_2$ is generically finite to its image, and assume both $X \to B_i$ surjective. 
\begin{enumerate}
\item Assume $B_1,B_2$ are of general type. Use Viehweg's theorem and the Kodaira defect inequality to conclude that $X$ is of general type. (Hint for a key step: construct a subvariety of general type $V \subset B_1$, such that $X \times_{B_1} V \to B_2$ is generically finite and surjective.) 
\item Assume $\kappa(B_1),\kappa(B_2) \geq 0$. Show that if Iitaka's conjecture holds true, then $\kappa(X) \geq 0$.
\end{enumerate}
\end{exercise}

\begin{exercise}\label{Ex:Lang-map}
Let $X$ be a smooth projective variety. Using the previous exercise, show that there is a dominant rational map 
$$L_X: X \das L(X)$$ 
such that 
\begin{enumerate}
\item $L(X)$ is of general type, and   
\item the map is universal: if $g: X \das Z$ is a dominant rational map with $Z$ of general type, there is a unique rational map $L(g):L(X) \das Z$ such that the following diagram commutes:
$$\xymatrix{
X \ar@{-->}^{L_X}[r]\ar@{-->}_g[rd] & L(X)\ar@{-->}[d]^{L(g)} \\ & Z.
}$$
\end{enumerate}
\end{exercise}

I call the map $L_X$  {\em the Lang map of $X$}, and $L(X)$ the {\em Lang variety of $X$}. 

\subsection{Uniruled varieties and rationally connected fibrations}
\subsubsection{Uniruled varieties} For simplicity let us assume here that $k$ is algebraically closed.

As indicated above, a variety $X$ is said to be {\em uniruled} if there is a $d-1$-dimensional variety $B$ and a dominant rational map $B\times \PP^1 \das X$. Instead of $B\times \PP^1$ one can take any variety $Y \to B$ whose generic fiber is a curve of genus 0. As discussed above, if $X$ is uniruled then $\kappa(X) = -\infty$. The converse is the important $(-\infty)$-Conjecture \ref{Conj:k<0-uniruled}.

A natural question is, can one ``take all these rational  curves out of the picture?" The answer is yes, in the best possible sense. 

\begin{definition} A smooth projective variety $P$ is said to be {\em rationally connected} if through any two points $x,y\in P$ there is a morphism from a rational curve $C \to P$ having $x$ and $y$ in its image.
\end{definition}
There are various equivalent ways to characterize rationally connected varieties. 
\begin{theorem}[Campana \cite{Campana-RC}, Koll\'ar-Miyaoka-Mori \cite{Kollar-Miyaoka-Mori}]
Let $P$ be a smooth projective variety. The following are equivalent:
\begin{enumerate}
\item $P$ is rationally connected.
\item Any two points are connected by a chain of rational curves.
\item For any finite set of points $S \subset P$, there is a morphism from a rational curve $C \to P$ having $S$ in its image.
\item There is a ``very free" rational curve on $P$ - if $\dim P>2$ this means there is a rational curve $C\subset P$ such that the normal bundle $N_{C\subset P}$ is ample.  
\end{enumerate}
\end{theorem}

Key properties:
\begin{theorem}[\cite{Campana-RC,Kollar-Miyaoka-Mori}] Let $X$ and $X'$ be smooth projective varieties, with $X$ rationally connected.
\begin{enumerate}
\item If  $X \das X'$ is a dominant rational map (in particular when $X$ and $X'$ are birationally equivalent) then $X'$ is rationally connected.
\item  If  $X'$ is deformation-equivalent to $X$ then $X'$ is rationally connected.
\item If  $X' = X_{k'}$ where $k'/k$ is an algebraically closed field extension,  then $X'$ is rationally connected if and only if $X$ is.
\end{enumerate}
\end{theorem}

\begin{exercise}
A variety is unirational if it is a dominant image of $\PP^n$. Show that every unirational variety is rationally connected.
\end{exercise}

On the other hand, one expects the following:
\begin{conjecture}[Koll\'ar]
There is a rationally connected threefold which is not unirational.  There should also exist some hypersurface of degree $n$ in $\PP^n$, $n\geq 4$ which is not unirational.
\end{conjecture}

Rational connectedness often arises when there is some negativity of differential forms, as in the following statement.  A smooth projective variety $X$ is Fano if its anti-canonical divisor is ample. We have the following:
\begin{theorem}[Koll\'ar-Miyaoka-Mori, Campana]
A Fano variety is rationally connected.
\end{theorem}

\begin{conjecture}[Koll\'ar-Miyaoka-Mori, Campana]
\begin{enumerate}
\item
A variety $X$ is rationally connected if and only if $$H^0(X, (\Omega^1_X)^{\otimes n}) = 0$$ for every positive integer $n$.
\item A variety $X$ is rationally connected if and only if every positive dimensional dominant image $X \das Z$ has $\kappa(Z) = -\infty$. 
\end{enumerate}
\end{conjecture}
This conjecture  follows from the minimal model program.

Now we can break any variety $X$ into a rationally connected fiber over a non-uniruled base:
\begin{theorem}[Campana, Koll\'ar-Miyaoka-Mori, Graber-Harris-Starr] Let $X$ be a smooth projective variety. There is a birational morphism $X' \to X$, a variety $Z(X)$, and a dominant morphism  $X' \to Z(X)$ with connected fibers,  such that 
\begin{enumerate}
\item The general fiber of $X' \to Z(X)$ is rationally connected, and 
\item $Z(X)$ is not uniruled. 
\end{enumerate}
Moreover, $X' \to X$ is an isomorphism in a neighborhood of the general fiber  of  $X' \to Z(X)$.
\end{theorem}
The existence of a fibration containing ``most"  rational curves was proven in the original papers by Campana and Koll\'ar-Miyaoka-Mori. The crucial fact that $Z(X)$ is not uniruled was proven by Graber, Harris and Starr in \cite{Graber-Harris-Starr}.

\subsubsection{} The rational map  $r_X:X \das Z(X)$ is called the {\em maximally rationally connected fibration of $X$} (or MRC fibration of $X$) and $Z(X)$, which is well defined up to birational equivalence, is called the {\em MRC quotient} of $X$. 

\subsubsection{} \label{Sec:other-fibrations} The MRC fibration has the universal property of being ``final" for dominant rational maps $X \to B$ with rationally connected fibers.

One can construct similar fibrations with similar universal property for maps with fibers having $H^0(X_b, (\Omega^1_{X_b})^{\otimes n}) = 0$, or for fibers having no dominant morphism to positive dimensional varieties of nonnegative Kodaira dimension. Conjecturally these agree with $r_X$. Also conjecturally, assuming Iitaka's conjecture,  there exists  $X \das Z'$ which is initial for maps to varieties of non-negative Kodaira dimension. This conjecturally will also  agree with $r_X$. All these conjecture would follow from the ``good minimal model" conjecture.

\subsubsection{Arithmetic, finally}
The set of rational points on a rational curve is Zariski dense. The following is a natural extension:

\begin{conjecture}[Campana] \label{Conj:rc-dense}
Let $P$ be a rationally connected variety over a number field $k$. Then rational points on $P$ are potentially dense.
\end{conjecture}

This conjecture and its sister below \ref{Conj:k=0-dense} was implicit in works of many, including Bogomolov, Colliot-Th\'el\`ene, Harris, Hassett,  Tschinkel. 

\subsection{Geometry and arithmetic of the Iitaka fibration}
We now want to understand the geometry and arithmetic of varieties such as $Z(X)$, i.e. non-uniruled varieties.  In view of Conjecture \ref{Conj:k<0-uniruled}, we focus on the case $\kappa(X) \geq 0$.

So let $X$ satisfy $\kappa(X) \geq 0$, and consider the Iitaka fibration $X \das I(X)$. The next proposition follows from easy subadditivity and Siu's theorem:
\begin{proposition}
Let $F$ be a general fiber of $X \to I(X)$. Then $\kappa(F) = 0$.
\end{proposition}

\begin{conjecture}[Campana]\label{Conj:k=0-dense}
Let $F$ be a variety over a number field $k$ satisfying $\kappa(F) = 0$. Then rational points on $F$ are potentially dense.
\end{conjecture}

\begin{exercise}
Recall the Lang map in \ref{Ex:Lang-map}. Assuming Conjecture \ref{Conj:k<0-uniruled}, show that $L(X)$ is the result of applying MRC fibrations and Iitaka fibrations until the result stabilizes.  
\end{exercise}

\subsection{Lang's conjecture}

In this section we let $k$ be a number field, or any field which is finitely generated over $\QQ$.

A highly inspiring conjecture in diophantine geometry is the following:

\begin{conjecture*}[Lang's conjecture, weak form]
Let $X$ be a smooth projective variety of general type over $k$. Then $X(k)$ is not Zariski-dense in $X$. 
\end{conjecture*}

In fact,  motivated by analogy with conjectures on the Kobayashi pseudo-metric  of a variety of general type, Lang even proposed the following:

\begin{conjecture*}[Lang's geometric conjecture]
Let $X$ be a smooth complex projective variety of general type. There is a Zariski closed proper subset $S(X)\subset X$, whose irreducible components are not of general type, and such that every irreducible subset $T\subset X$ not of general type is contained in $S(X)$.
\end{conjecture*}

The notation ``$S(X)$" stands for ``the special subvariety of $X$". It is not hard to see that $S(X)$ is defined over any field of definition of $X$.
The two conjectures combine to give:

\begin{conjecture*}[Lang's conjecture, strong form]
Let $X$ be a smooth projective variety of general type over  $k$. Then for any finite extension $k'/k$, the set  $(X\setmin S(X))(k')$ is finite. 
\end{conjecture*}

Here is a simple consequence: 
\begin{proposition}
Assume Lang's conjecture holds true. Let $X$ be a smooth projective variety over a number field $k$. Assume there is a dominant rational map $X \to Z$, such that $Z$ is a positive dimensional variety of general type (i.e.,  $\dim L(X)>0$). Then $X(k)$ is not Zariski-dense  in $X$.  
\end{proposition}

\subsection{Uniformity of rational points.} Lang's conjecture can be investigated whenever one has a variety of general type around. By considering certain subvarieties of the moduli space $\cM_{g,n}$ of curves of genus $g$ with $n$ distinct points on them, rather surprising and inspiring implications on the arithmetic of curves arise. This is the subject of the work \cite{Caporaso-Harris-Mazur} of L. Caporaso, J. Harris and B. Mazur. Here are their key results:

\begin{theorem}
Assume that the weak Lang's conjecture holds true. Let $k$ be as above, and let $g>1$ be an integer. Then there exists an integer $N(k,g)$ such that for every algebraic curve $C$ of genus $g$ over $k$ we have $$\# C(k) \leq   N(k,g).$$
\end{theorem}

\begin{theorem}
Assume that the strong Lang's conjecture holds true. Let  $g>1$ be an integer. Then there exists an integer $N(g)$ such that for every finitely generated field $k$ there are, up to isomorphisms, only finitely many algebraic curves $C$ of genus $g$ over $k$  with $\# C(k) >   N(g).$
\end{theorem}

Further results along these lines, involving higher dimensional varieties and involving stronger results on curves can be found in \cite{Hassett}, \cite{Abramovich-quadratic}, \cite{Pacelli}, \cite{Abramovich-Voloch}, \cite{Abramovich-fibered}. For instance, P. Pacelli's result in \cite{Pacelli} says that the number $N(k,g)$ can be replaced for number fields by $N(d,g)$, where $d=[k:\QQ]$.

The reader may decide whether this shows the great power of the conjectures or their unlikelihood. I prefer to be agnostic and rely on the conjectures for inspiration.

\subsection{The search for an arithmetic dichotomy} 

As demonstrated in table  \ref{Tab:rational}, potential density of rational points on curves is dictated by geometry. Lang's conjecture carves out a class of higher dimensional varieties for which rational points are, conjecturally, not potentially dense. Can this be extended to a dichotomy as we have for curves?

One can naturally wonder - is the Kodaira dimension itself enough for determining potential density of points? Or else, maybe just the inexistence of a map to a positive dimensional variety of general type?

\subsubsection{Rational points on surfaces} The following table, which I copied from a lecture of L. Caporaso \cite{Caporaso}, describes what is known about surfaces. 

\begin{center} {\sc Caporaso's table: rational points on surfaces} \\[5mm]   \begin{tabular}{|l|l|l|}\hline 
 Kodaira dimension &  $X(k)$ potentially dense &  $X(k)$ never dense \\ \hline\hline
$\kappa= -\infty$ & $\PP^2$ & $\PP^1\times C$ {\small($g(C)\geq 2$)} \\ \hline 
 $\kappa=0$ & $E\times E$, many others & {\bf none known} \\ \hline
 $\kappa=1$ &  many examples & $E \times C$ {\small($g(C)\geq 2$)} \\ \hline
 $\kappa=2$ & {\bf none known} & many examples \\ \hline
\end{tabular} \end{center}

The bottom row is the subject of Lang's conjecture, and the $\kappa=0$ row is the subject of Conjecture \ref{Conj:k=0-dense}.

\subsubsection{Failure of the dichotomy using $\kappa(X)$}
 
 The first clear lesson we learn from this is, as Caporaso aptly put it in her lecture, that 
 \begin{center} \begin{tabular}{|c|}\hline 
  diophantine geometry is not governed by the Kodaira dimension. \\ \hline
\end{tabular} \end{center}
 On the top row we see that clearly: on a ruled surface over a curve of genus $\geq 2$, rational points can never be dense by Faltings's theorem. So it behaves very differently from a rational surface. 
 
Even if one insists on working with varieties of non-negative Kodaira dimension, the $\kappa=1$ row gives us trouble. 
\begin{exercise*}\label{Ex:special-elliptic}

 Take a Lefschetz pencil of cubic curves in $\PP^2$, parametrized by $t$, and assume that it has two sections $s_1, s_2$ whose difference is not torsion on the generic fiber. We use $s_1$ as the origin.
\begin{enumerate}
\item Show that the dualizing sheaf of the total space $S$ is $\cO_S(-[F])$, where $F$ is a fiber. 
\item Show that the relative dualizing sheaf is $\cO_S([F])$.
Take the base change $t = s^3$. We still have two sections, still denoted $s_1,s_2$, such that the difference is not torsion. We view $s_1$ as origin.

Show that the relative dualizing sheaf of the new surface $X$ is  $\cO_X(3[F])$ and its dualizing sheaf is $\cO_X([F])$.
 Conclude that 
the resulting  surface $X$ has Kodaira dimension 1.
\item For any rational point $p$ on $\PP^1$ where the section $s_2$ of $X \to \PP^1$ is not torsion, the fiber has a dense set of rational points. 

In characteristic 0 it can be shown that the set of such points is dense. For instance, by Mazur's theorem the rational torsion points have order at most 12, and therefore they lie on finitely many points of intersection of $s_2$ with the locus of torsion points of order $\leq 12$.
\item Conclude that $X$ has a dense set of rational points.
 \end{enumerate}
 \end{exercise*}

\subsubsection{Failure of the dichotomy using the Lang map}

The examples given above still allow for a possible dichotomy based on the existence of a nontrivial map to a variety of general type. But the following example, which fits on the right column on row $\kappa=1$, shows this doesn't work either. The example is due to Colliot-Th\'el\`ene, Skorobogatov and Swinnerton-Dyer \cite{CT-S-SD}.

\begin{example*}\label{Ex:non-special-elliptic}
Let $C$ be a curve with an involution $\phi: C\to C$, such that the quotient  is rational. Consider an elliptic curve $E$ with a 2-torsion point $a$, and consider the fixed-point free action of $\ZZ/2\ZZ$ on $Y = E \times C$ given by: 
$$(x,y) \mapsto (x+a, \phi(y)).$$

Let the quotient of $Y$ by the involution be $X$. Then $L(X)$ is trivial, though rational points on $X$ are not potentially dense by Chevalley-Weil and Faltings. 
\end{example*}

In the next lecture we address a conjectural approach to a dichotomy - due to F. Campana - which has a chance to work.

\subsection{Logarithmic Kodaira dimension and the Lang-Vojta conjectures}

We now briefly turn our attention to open varieties, following the lesson in section \ref{Lesson:open}. 

Let $\oX$ be a smooth projective variety, $D$ a reduced normal crossings divisor. We can consider the quasiprojective variety $X = \oX \setminus D$.

The logarithmic Kodaira dimension of $X$ is defined to be  the Iitaka dimension $\kappa(X):=\kappa(\oX, K_\oX+D)$. We say that $X$ is of {\em logarithmic general type} if $\kappa(X)= \dim X$.

It can be easily shown that $\kappa(X)$ is independent of the completion $X \subset \oX$, as long as $\oX$ is smooth and $D$ is a normal crossings divisor. More invariance properties can be discussed, but will take us too far afield.

Now to arithmetic: suppose $\cX$ is a model of $X$ over $\OKS$.
We can consider integral points $\cX(\cO_{L,S_L})$ for any finite extension $L/k$ and enlargement $S_L$ of the set of places over $S$. 

The Lang-Vojta conjecture is the following:
\begin{conjecture} If $X$ is of logarithmic general type, then integral points are not potentially dense on $X$, i.e. $\cX(\cO_{L,S_L})$ is not Zariski dense for any $L,S_L$.
\end{conjecture}

\subsubsection{} In case $X = \oX$ is already projective, the Lang-Vojta conjecture reduces to Lang's conjecture: $X$ is simply a variety of general type, integral points on $X$ are the same as rational points, and Lang's conjecture asserts that $X(k)$ is not Zariski-dense in $X$. 

\subsubsection{} The Lang-Vojta conjecture turns out to be a particular case of  a more precise and more refined conjecture of Vojta, which will be discussed in a later lecture.

\section{Campana's program}

For this section one important road sign is
\begin{center}\fbox{\parbox{3.2in}{\large THIS SITE IS UNDER CONSTRUCTION\\ DANGER! HEAVY EQUIPMENT CROSSING}}
\end{center}

A quick search on the web shows close to the top a number of web sites deriding the idea of ``site under construction". Evidently these people have never engaged in research!

\subsubsection{} Campana's program is a new method of breaking algebraic varieties into ``pieces" which builds upon Iitaka's program, but, by using a particular structure on varieties which I will call ``Campana constellations" enables one to get closer to a classification which is compatible with arithmetic properties. There is in fact an underlying more refined structure which I call ``firmament" for the Campana constellation, which might be the more fundamental structure to study. It truly does say something about rational points.

\subsubsection{} The term ``constellation" is inspired by Aluffi's celestial \cite{Aluffi}, which is in turn inspired by Hironaka. 

Campana used the term ``orbifold", in analogy to orbifolds  used in geometry, but the analogy breaks very early on. A suggested replacement ``orbifold pair" still does not make me too happy. Also, ``Campana pair" is a term which Campana himself is not comfortable using, nor could he shorten it to just ``pair", which is insufficient. I was told by Campana that he would be happy to use ``constellations" if the term catches.

\subsection{One dimensional  Campana constellations}

\subsubsection{The two key examples: elliptic surfaces}
Let us inspect again Caporaso's table of surfaces, and concentrate on $\kappa=1$. We have in \ref{Ex:special-elliptic} and \ref{Ex:non-special-elliptic} two examples, say $S_1 \to \PP^1$ and  $S_2 \to \PP^1$ of elliptic surfaces of Kodaira dimension 1 fibered over $\PP^1$. But their arithmetic behavior is very different. 

Campana asked the question: is there an underlying structure on the base $\PP^1$ from which we can deduce this difference of behavior?

The key point is that the example in \ref{Ex:non-special-elliptic} has $2g+2$ double fibers lying over a divisor $D\subset \PP^1$. This means that the elliptic surface $S_2\to \PP^1$ can be lifted to $S_2 \to \cP$, where $\cP$ is the orbifold structure  $\PP^1(\sqrt{D})$ on $\PP^1$ obtained by taking the square root of $D$. Following the ideas of Darmon and Granville in \cite{Darmon-Granville}, one should consider the canonical divisor class $K_\cP$ of $\cP$, viewed as a divisor with rational coefficients on $\PP^1$, namely $K_{\PP^1} + (1-1/2) D$. In general, when one has an $m$-fold fiber over a divisor $D$, one wants to take $D$ with coefficient $(1-1/m)$.

Darmon and Granville prove, using Chevalley-Weil and Faltings, that such an orbifold $\cP$ has potentially dense set of integral points if and only if the Kodaira dimension $\kappa(\cP) = \kappa(\cP, K_\cP) < 1$. And the image of a rational point on $S_2$ is an integral point on $\cP$. This fully explains our example: since integral points on $\cP = \PP^1(\sqrt{D})$ are not Zariski dense, and since rational points on $S_2$ map to integral points on $\cP$, rational points on $S_2$ are not dense.

\subsubsection{The multiplicity divisor} What should we declare the structure to be when we have a fiber that looks like $x^2y^3=0$, i.e. has two components of multiplicities 2 and 3? Here Campana departs from the classical orbifold picture: the highest classical orbifold to which the fibration lifts has no new structure lying under such a fiber, because $\gcd(2,3)=1$. Campana makes a key observation that a rich and interesting classification theory arises if one instead considers $\min(2,3) = 2$ as the basis of the structure. 

\begin{definition}[Campana]
Consider a dominant morphism $f:X \to Y$ with $X, Y$ smooth and $\dim Y = 1$. Define a divisor with rational coefficients $\Delta_f = \sum \delta_p p$ on $Y$ as follows:  assume the divisor $f^{*} p$ on $X$ decomposes as  $f^{*} p = \sum m_iC_i$, where $C_i$ are the distinct irreducible components of the fiber taken with reduced structure. Then set $$\delta_p \ = \ 1-\frac{1}{m_p}, \ \ \  \text{ where  }\ \ \    m_p\  =\  \min_i\, m_i. $$
\end{definition}

\begin{definition}[Campana] 
\begin{enumerate}
\item A  Campana constellation curve  $(Y/\Delta)$ is a pair consisting of a curve $Y$ along with a divisor $\Delta= \sum \delta_p p$ with rational coefficients, where each $\delta_p$ is of the form $\delta_p  =  1-1/m_p$ for some integer $m_p$.

\item The  Campana constellation base of $f: X \to Y$ is the structure pair consisting of $Y$ with the divisor $\Delta_f$ defined above, denoted $(Y/\Delta_f)$.
\end{enumerate}
\end{definition}

The word used by Campana is {\em orbifold}, but as I have argued,  the analogy with orbifolds is  shattered  in this very definition. 


The new terminology ``constellation" will become better justified and much more laden with meaning when we consider $Y$ of higher dimension.

Campana's definition deliberately does not distinguish between the structure coming from a fiber of type $x^2=0$ and one of type $x^2y^3=0$. We will see later a way to resurrect the difference to some extent using the notion of {\em firmament}, by which a Campana constellation hangs.

\begin{definition}[Campana] The Kodaira dimension of a  Campana constellation curve $(Y/\Delta)$ is defined as the following Iitaka dimension:
$$\kappa\left(\,(Y/\Delta)\,\right) = \kappa(Y, K_Y + \Delta).$$ We say that $(Y/\Delta)$ is of general type if it has Kodaira dimension 1. We say that it is {\em special} if it is not of general type.
\end{definition}

\begin{exercise}
Classify  {\em special}  Campana constellation curves over $\CC$. See \cite{Campana2} for a detailed discussion.
\end{exercise}

\subsubsection{Models and integral points} Now to arithmetic. As we learned in Lesson \ref{Lesson:open}, when dealing with a variety with a structure given by a divisor, we need to speak about {\em integral points} on an {\em integral model} of the structure. Thus let $\cY$ be an integral model of $Y$, proper over $\OKS$, and denote by $\widetilde\Delta$ the closure of $\Delta$.  As above we assume that $\Delta$ is the union of integral points, which we denote now $z_i$ of $Y$, and for simplicity let us assume that they are disjoint (we can always achieve  this by enlarging $S$). It turns out that there is more than one natural notion to consider in our theory - soft and firm. The firm notion will be introduced when firmaments are considered.

\begin{definition*} A $k$ rational point $x$ on $Y$, considered as an
  integral point  of $\cY$, is said to be a {\em soft
    $S$-integral points on $(\cY/\widetilde\Delta)$} if for any integral point $z$ in $\widetilde \Delta$, and any
  nonzero prime $\wp\subset \OKS$ such that  the reductions coincide: $x_\wp=z_\wp$, we have $$\operatorname{mult}_\wp(
  x\cap z) \geq  m_p.$$
\end{definition*}

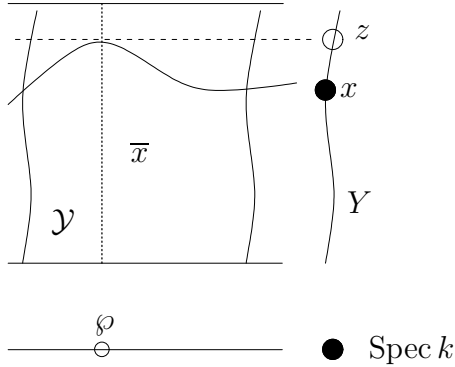
\begin{center}
\setlength{\unitlength}{0.0005in}
\begingroup\makeatletter\ifx\SetFigFont\undefined%
\gdef\SetFigFont#1#2#3#4#5{%
  \reset@font\fontsize{#1}{#2pt}%
  \fontfamily{#3}\fontseries{#4}\fontshape{#5}%
  \selectfont}%
\fi\endgroup%
{\renewcommand{\dashlinestretch}{30}
\begin{picture}(4285,3906)(0,-10)
\put(3387,3504){\ellipse{212}{212}}
\put(3312,2979){\blacken\ellipse{212}{212}}
\put(3312,2979){\ellipse{212}{212}}
\put(3387,279){\blacken\ellipse{212}{212}}
\put(3387,279){\ellipse{212}{212}}
\put(987,279){\ellipse{150}{150}}
\path(12,3879)(2862,3879)
\drawline(2862,3879)(2862,3879)
\path(12,1179)(2862,1179)
\path(12,279)(2862,279)
\dottedline{45}(987,3879)(987,1179)
\path(3462,3804)(3462,3802)(3460,3797)
	(3458,3787)(3455,3773)(3451,3753)
	(3446,3727)(3439,3695)(3431,3658)
	(3423,3617)(3414,3573)(3405,3526)
	(3395,3478)(3386,3429)(3377,3381)
	(3368,3334)(3360,3288)(3353,3244)
	(3346,3202)(3339,3162)(3334,3124)
	(3329,3088)(3325,3054)(3321,3021)
	(3318,2989)(3316,2958)(3314,2928)
	(3313,2899)(3312,2870)(3312,2841)
	(3312,2811)(3313,2780)(3314,2749)
	(3316,2718)(3318,2686)(3320,2654)
	(3323,2621)(3327,2588)(3330,2555)
	(3334,2521)(3338,2487)(3342,2453)
	(3347,2418)(3351,2384)(3356,2350)
	(3360,2316)(3365,2283)(3369,2250)
	(3374,2217)(3378,2185)(3382,2154)
	(3385,2124)(3388,2094)(3391,2065)
	(3394,2036)(3396,2009)(3398,1981)
	(3400,1954)(3401,1925)(3401,1896)
	(3401,1867)(3400,1837)(3399,1807)
	(3397,1776)(3395,1743)(3392,1709)
	(3388,1673)(3384,1635)(3378,1595)
	(3373,1553)(3367,1510)(3360,1466)
	(3353,1422)(3346,1379)(3339,1338)
	(3333,1300)(3327,1266)(3322,1237)
	(3318,1215)(3315,1198)(3314,1188)
	(3312,1182)(3312,1179)
\path(2637,3804)(2637,3802)(2635,3797)
	(2633,3787)(2630,3773)(2626,3753)
	(2621,3727)(2614,3695)(2606,3658)
	(2598,3617)(2589,3573)(2580,3526)
	(2570,3478)(2561,3429)(2552,3381)
	(2543,3334)(2535,3288)(2528,3244)
	(2521,3202)(2514,3162)(2509,3124)
	(2504,3088)(2500,3054)(2496,3021)
	(2493,2989)(2491,2958)(2489,2928)
	(2488,2899)(2487,2870)(2487,2841)
	(2487,2811)(2488,2780)(2489,2749)
	(2491,2718)(2493,2686)(2495,2654)
	(2498,2621)(2502,2588)(2505,2555)
	(2509,2521)(2513,2487)(2517,2453)
	(2522,2418)(2526,2384)(2531,2350)
	(2535,2316)(2540,2283)(2544,2250)
	(2549,2217)(2553,2185)(2557,2154)
	(2560,2124)(2563,2094)(2566,2065)
	(2569,2036)(2571,2009)(2573,1981)
	(2575,1954)(2576,1925)(2576,1896)
	(2576,1867)(2575,1837)(2574,1807)
	(2572,1776)(2570,1743)(2567,1709)
	(2563,1673)(2559,1635)(2553,1595)
	(2548,1553)(2542,1510)(2535,1466)
	(2528,1422)(2521,1379)(2514,1338)
	(2508,1300)(2502,1266)(2497,1237)
	(2493,1215)(2490,1198)(2489,1188)
	(2487,1182)(2487,1179)
\path(312,3804)(312,3802)(310,3797)
	(308,3787)(305,3773)(301,3753)
	(296,3727)(289,3695)(281,3658)
	(273,3617)(264,3573)(255,3526)
	(245,3478)(236,3429)(227,3381)
	(218,3334)(210,3288)(203,3244)
	(196,3202)(189,3162)(184,3124)
	(179,3088)(175,3054)(171,3021)
	(168,2989)(166,2958)(164,2928)
	(163,2899)(162,2870)(162,2841)
	(162,2811)(163,2780)(164,2749)
	(166,2718)(168,2686)(170,2654)
	(173,2621)(177,2588)(180,2555)
	(184,2521)(188,2487)(192,2453)
	(197,2418)(201,2384)(206,2350)
	(210,2316)(215,2283)(219,2250)
	(224,2217)(228,2185)(232,2154)
	(235,2124)(238,2094)(241,2065)
	(244,2036)(246,2009)(248,1981)
	(250,1954)(251,1925)(251,1896)
	(251,1867)(250,1837)(249,1807)
	(247,1776)(245,1743)(242,1709)
	(238,1673)(234,1635)(228,1595)
	(223,1553)(217,1510)(210,1466)
	(203,1422)(196,1379)(189,1338)
	(183,1300)(177,1266)(172,1237)
	(168,1215)(165,1198)(164,1188)
	(162,1182)(162,1179)
\dashline{60.000}(3162,3504)(87,3504)
\path(12,2829)(14,2831)(18,2835)
	(26,2842)(38,2853)(55,2869)
	(77,2889)(104,2914)(136,2942)
	(172,2975)(211,3010)(254,3047)
	(298,3085)(343,3123)(388,3161)
	(433,3198)(477,3233)(519,3266)
	(560,3298)(599,3326)(637,3352)
	(673,3376)(707,3396)(739,3415)
	(770,3431)(800,3444)(829,3455)
	(857,3464)(884,3471)(910,3475)
	(936,3478)(962,3479)(989,3478)
	(1017,3475)(1045,3471)(1072,3464)
	(1101,3456)(1129,3446)(1158,3435)
	(1187,3422)(1217,3408)(1248,3392)
	(1278,3376)(1310,3358)(1341,3339)
	(1373,3319)(1405,3299)(1437,3279)
	(1469,3258)(1501,3237)(1533,3216)
	(1565,3195)(1597,3175)(1628,3155)
	(1660,3136)(1691,3118)(1721,3100)
	(1752,3084)(1783,3069)(1813,3054)
	(1844,3041)(1875,3029)(1904,3019)
	(1933,3010)(1964,3002)(1995,2995)
	(2027,2989)(2060,2984)(2095,2980)
	(2132,2978)(2170,2976)(2211,2975)
	(2254,2976)(2298,2977)(2346,2979)
	(2395,2982)(2447,2986)(2500,2990)
	(2554,2996)(2609,3001)(2664,3007)
	(2717,3014)(2769,3020)(2817,3026)
	(2861,3032)(2899,3038)(2933,3042)
	(2960,3046)(2980,3049)(2995,3051)
	(3005,3053)(3010,3054)(3012,3054)
\put(3762,129){\makebox(0,0)[lb]{{\SetFigFont{12}{14.4}{\rmdefault}{\mddefault}{\updefault}$\Spec k$}}}
\put(3462,2904){\makebox(0,0)[lb]{{\SetFigFont{12}{14.4}{\rmdefault}{\mddefault}{\updefault}$x$}}}
\put(1287,2229){\makebox(0,0)[lb]{{\SetFigFont{12}{14.4}{\rmdefault}{\mddefault}{\updefault}$\overline
      x$}}}
\put(462,1479){\makebox(0,0)[lb]{{\SetFigFont{12}{14.4}{\rmdefault}{\mddefault}{\updefault}$\cY$}}}
\put(3612,3504){\makebox(0,0)[lb]{{\SetFigFont{12}{14.4}{\rmdefault}{\mddefault}{\updefault}$z$}}}
\put(912,429){\makebox(0,0)[lb]{{\SetFigFont{12}{14.4}{\rmdefault}{\mddefault}{\updefault}$\wp$}}}
\put(3537,1704){\makebox(0,0)[lb]{{\SetFigFont{12}{14.4}{\rmdefault}{\mddefault}{\updefault}$Y$}}}
\end{picture}
}
\end{center}

A key property of this definition is:

\begin{proposition} Assume $f:X \to Y$ extends to a good model  $\tilde f: \cX \to \cY$. Then
the image of a rational point on $X$ is a soft $S$-integral point on $(\cY/\widetilde\Delta_f)$.
\end{proposition}

So rational points on $X$ can be investigated using integral points on a model of $Y$. This makes the following very much relevant:

\begin{conjecture}[Campana] \label{Conj:Cc-curve}
Suppose the Campana constellation curve $(Y/\Delta)$ is of general type. Then the set of  soft $S$-integral point on any model $\cY$ is not Zariski dense.
\end{conjecture}

This conjecture is not likely to follow readily from Faltings's theorem, as the following example suggests. 

\begin{example} Let $n \geq 4$ be an integer. Let $Y\simeq \PP^1$ and $\Delta$ the divisor supported at $0, 1$ and $\infty$  with all multiplicities equal to $(n-1)/n$. Then $(Y/\Delta)$ is of general type.

Using the same as a model over $\Spec \ZZ$, we see that a point $y$ on $Y$ is a soft integral point on $(Y/\Delta)$ if at every prime where $y$ reduces to $0,1$ or $\infty$, the multiplicity of this reduction is at least $n$. 

Considering a triple $a,b,c$ of relatively prime integers with $a^N+b^N=c^N$, 
the point $(a^N:c^N)$ on $Y$ is a soft integral point as soon as $N\geq n$.

 It follows that Campana's conjecture \ref{Conj:Cc-curve} implies asymptotic Fermat over any number field. 
\end{example}

It also seems that the conjecture does not follow readily from any of the methods surrounding Wiles's proof of Fermat.
As we'll see in the last lecture, the conjecture does follow from the $abc$ conjecture (which implies asymptotic Fermat). In particular we have the following theorem in the function field case.

\begin{theorem}[Campana] Let $B$ be a complex algebraic curve, and $K$ its function field, and let $S \subset B$ be a finite set of closed points.
Let $(Y/\Delta)$ be a Campana constellation curve of general type defined over the function field $K$. Then the set of non-constant  soft $S$-integral points on any model $\cY \to B$ is not Zariski dense.
\end{theorem}

\subsubsection{Some examples} 
\begin{enumerate}
\item Consider $f: \AA^2 \to \AA^1$ given by $t = x^2$. The constellation base has $\Delta = 1/2(0)$, where $(0) $ is the origin on $\AA^1$. Sections of $\cO_Y(K+\Delta)$ are generated by $dt$, sections of $\cO_Y\left(2(K+\Delta)\right)$ by $(dt)^2/t$, and sections of $\cO_Y\left(3(K+\Delta)\right)$ by $(dt)^3/t$.
\item The same structure occurs for $f: \AA^2 \to \AA^1$ given by $t = x^2y^3$ and $f: \AA^2 \to \AA^1$ given by $t = x^2y^2$.
\item For $f: \AA^2 \to \AA^1$ given by $t = x^2y$ the constellation base is trivial. 
\item for $f: \AA^2 \to \AA^1$ given by $t = x^3y^4$, we get $\Delta = 2/3(0)$. Again sections of $\cO_Y\left(2(K+\Delta)\right)$ are generated by $(dt)^2/t$, but sections of $\cO_Y\left(3(K+\Delta)\right)$ are generated by $(dt)^3/t^2$.
\end{enumerate}

\subsection{Higher dimensional Campana constellations} We turn now to the analogous situation of $f:X \to Y$ with higher dimensional $Y$.

One seeks to define objects, say Campana constellations $(Y/\bDelta)$, in analogy to the case of curves, which in some sense should help us understand the geometry and arithmetic of plain varieties mapping to them. 

Ideally, these objects should form a category extending the category of varieties, at least with some interesting class of morphisms. Ideally these objects should have a good notion of differential forms which fits into the standard theory of birational geometry, for instance having well-behaved Kodaira dimension. Ideally there should be a notion of integral points on $(Y/\bDelta)$ which says something about rational points of a ``plain" variety  $X$ whenever $X$ maps to $(Y/\bDelta)$.

The theory we describe in this section, which is due to Campana in all but some details, relies on divisorial data. We will describe a category of objects, called Campana constellations, which at the moment only allows dominant morphisms. This means that we do not have a satisfactory description of integral points, since integral points are sections, and sections are not dominant morphisms. The  theory of firmaments aims at resolving this problem.

 Unfortunately, points on $Y$ are no longer divisors. And divisors on $Y$ are not quite sufficient to describe codimension $>1$ behavior. Campana resolves this by considering all birational models of $Y$ separately. This brings him to define various invariants, such as Kodaira dimension, depending on a morphism $X \to Y$ rather than of the structure $(Y/\bDelta)$ itself. I prefer to put all this data together using the notion of a b-divisor, introduced by Shokurov \cite{Shokurov}, based on ideas by Zariski \cite{Zariski}.  I was also inspired by Aluffi's  \cite{Aluffi}.
The main advantage is that all invariants will be defined directly on the level of $(Y/\bDelta)$. This structure has the disadvantage that it is not obviously computable in finite or combinatorial terms. It turns out that it is - this again will be addressed using firmaments.

\begin{definition} Let $k$ be a field and $Y$ a variety over $k$. A {\em rank 1 discrete valuation} on the function field $\cK =\cK(Y)$ over $k$ is a surjective group homomorphism $\nu:\cK^\times \to \ZZ$, sending $k^\times$ to 0, satisfying 
$$\nu(x+y) \geq \min(\nu(x),\nu(y))$$
 with equality unless $\nu(x) = \nu(y)$. We define $\nu(0) =\, +\infty$.
 
 The {\em valuation ring of $\nu$} is defined as $$R_\nu \ \ =\ \  \left\{x\in \cK\  \big|\  \nu(x) \geq 0\right\}.$$ Denote by $Y_\nu = \Spec R_\nu$, and its unique closed point $s_\nu$.
 
  A rank 1 discrete valuation $\nu$ is {\em divisorial} if there is a birational model $Y'$ of $Y$ and an irreducible  divisor $D'\subset Y'$ such that for all nonzero $x \in \cK(Y) = \cK(Y')$ we have  $$\nu(x) \ = \ \operatorname{mult}_{D'}x.$$ In this case we say $\nu$ has divisorial center $D'$ in $Y'$.
\end{definition}

\begin{definition}
A b-divisor $\bDelta$ on $Y$ is an expression of the form 
$$\bDelta\ \  =\ \  \sum_\nu\, c_\nu\cdot \nu,$$ a possibly infinite sum over divisorial valuations of $\cK(Y)$ with rational coefficients, which satisfies the following finiteness condition:
\begin{itemize}\item  for each birational model $Y'$ there are only finitely many $\nu$ with divisorial center on $Y'$ having $c_\nu \neq 0$. \end{itemize}

A b-divisor is of {\em orbifold type} if for each $\nu$ there is a positive integer $m_\nu$ such that $c_\nu= 1-1/m_\nu$. 
\end{definition}

Before we continue, here is an analogue of the strict transform of a divisor: 
\begin{definition} Let $Y$ be a variety, $X$ a reduced scheme, and   let $f: X \to Y$ be a  morphism. Consider an integral scheme  $Y'$ with generic point $\eta$, a dominant morphism $Y' \to Y$, and the pullback $X \times_Y Y' \to Y'$. The {\em main part} $\widetilde{X \times_Y Y'}$ of  $X \times_Y Y'$ is the closure of the generic fiber $X \times_Y \eta$ inside  $X \times_Y Y'$. 
\end{definition}

Here is a higher dimensional  analogue of the divisor underlying the constellation curve of a morphism $X \to Y$:

\begin{definition}\label{Def:C-b-divisor} Let $Y$ be a variety, $X$ a reduced scheme, and   let $f: X \to Y$ be a  morphism, surjective on each irreducible component of $X$. For each divisorial valuation $\nu$ on $\cK(Y)$ consider $f': X_\nu' \to Y_\nu$, where $X'_\nu$ is a desingularization of the main part of the pullback $X \times_YY_\nu$.  Write ${f'}^*s_\nu = \sum m_i C_i$.  Define 
$$\delta_\nu \,=\, 1-\frac{1}{m_\nu} \quad \text{with} \quad m_\nu = \min_i m_i.$$ 

The {\em Campana b-divisor} on $Y$ associated to a dominant map $f: X \to Y$ is defined to be the b-divisor  $$\bDelta_f=\sum \delta_\nu \nu .$$
\end{definition}

\begin{exercise}
The definition is independent of the choice of desingularization $X'_\nu$.
\end{exercise}

This makes the b-divisor $\bDelta_f$ invariant under proper birational transformations on $X$ and $Y$. In particular the notion makes sense for a dominant rational map $f$.

\begin{definition}
\begin{enumerate}
\item A  {\em Campana constellation} $(Y/\bDelta)$ consists of  a variety $Y$ with a b-divisor $\bDelta$ such that, locally in the \'etale topology on $Y$, there is $f: X \to Y$ with $\bDelta=\bDelta_f$. 
\item The Campana constellation base of a morphism $X \to Y$ as above is $(Y/\bDelta_f)$ 
\item The trivial constellation on $Y$ is given by the zero b-divisor.
\item For each birational model $Y'$, define the $Y'$-divisorial part of $\bDelta$: 
$$\Delta_{Y'} = \sum_{\nu \text{ with divisorial support on } Y'} \delta_\nu \nu.$$
\end{enumerate} 
\end{definition}

The definition of a constellation feels a bit unsatisfactory because it requires, at least locally,  the existence of a morphism $f$. But using the notion of firmament, especially toroidal firmament, we will make this structure more combinatorial, in such a way that the existence of  $f$ is automatic.

\subsubsection{} Here's why I like the word ``constellation": think of a divisorial valuation $\nu$ as a sort of ``generalized point" on $Y$. Putting  $\delta_\nu>0$ suggests viewing a ``star" at that point. Replacing $Y$ by higher and higher models $Y'$ is analogous to using stronger and stronger telescopes to view farther stars deeper into space. The picture I have in my mind is somewhat reminiscent of the astrological meaning of ``constellation", not as just one group of stars, but rather as the arrangement of the entire heavens at the time the ``baby" $X \to Y$ is born. But hopefully it is better grounded in reality.

We now consider morphisms. For constellations we work only with dominant morphisms.
\begin{definition}
\begin{enumerate}
\item Let $(X/\bDelta_X)$ be a Campana constellation, and $f: X \to Y$ a {\em proper} dominant morphism. The constellation base $(Y, \bDelta_{f,\bDelta_X})$ is defined as follows: for each divisorial valuation $\nu$ of $Y$ and each divisorial valuation $\mu$ of $X$ with center $D$ dominating the center $E$ of $\nu$, let $$m_{\mu/\nu} = m_\mu \cdot \operatorname{mult}_D(f^*E).$$  Define
 $$m_\nu = \min_{\mu/\nu} m_{\mu/\nu}\  \text{ and }\ \delta_\nu = 1-\frac{1}{m_\nu}.$$ Then set as before
 $$\bDelta_{f,\bDelta_X} = \sum_\nu \delta_\nu \nu.$$
 \item Let $(X/\bDelta_X)$ and $(Y/\bDelta_Y)$ be Campana constellations and $f: X \to Y$ a dominant  morphism. Then $f$ is said to be a {\em constellation morphism} if for every divisorial valuation $\nu$ on $Y$ and any $\mu/\nu$ we have $m_\nu \leq m_{\mu/\nu}$, where as above $m_{\mu/\nu} = m_\mu \cdot \operatorname{mult}_D(f^*E)$. When $f$ is proper this just means  $\bDelta_Y \leq    \bDelta_{f,\bDelta_X}$.
\end{enumerate}
\end{definition}

Now to differential forms:

\begin{definition} A rational $m$-canonical differential $\omega$ on $Y$ is said to be regular on $(Y/\bDelta)$ if for every divisorial valuation $\nu$ on $\cK(Y)$, the polar multiplicity of $\omega$ at $\nu$ satisfies $$(\omega)_{\infty, \nu} \leq m \delta_\nu.$$ In other words, $\omega$ is a section of $\cO_{Y'}(m(K_{Y'} +\Delta_{Y'}))$ on every birational model $Y'$.  

The Kodaira dimension $\kappa(\,(Y/\bDelta)\,)$ is defined using the ring of regular $m$-canonical differentials on $(Y/\bDelta)$.
\end{definition}

\begin{exercise}
This is a birational invariant: if $Y$ and $Y'$ are proper and have the same function field, then $\kappa(\,(Y/\bDelta)\,) = \kappa(\,(Y'/\bDelta)\,)$.
\end{exercise}

\begin{theorem}[Campana \cite{Campana} Section 1.3]\label{Th:Campana-admissible}
There is a birational model $Y'$ with $\Delta_{Y'}$  a normal crossings divisor such that  $$\kappa(\,(Y/\bDelta)\,) = \kappa(Y', K_{Y'}+\Delta_{Y'}),$$ and moreover the algebra of regular pluricanonical differentials on $(Y/\bDelta)$ agrees with the algebra of sections $\oplus_{m\geq 0} H^0(Y', \cO_{Y'}(m(K_{Y'}+\Delta_{Y'})))$.
\end{theorem}

Campana calls such a model {\em admissible}. This is proven using Bogomolov sheaves, an important notion which is a bit far afield for the present discussion. The formalism of firmaments, especially toroidal firmaments, allows one to give a combinatorial proof of this result.

We remark that this theorem means that the new and ground-breaking finite generation theorem of \cite{BCHM} applies, so the algebra of regular pluricanonical differentials on $(Y/\bDelta)$ is finitely generated.

It is not difficult to see that any birational  model lying over an admissible model is also admissible.

\begin{definition}
A   Campana constellation  $(Y/\bDelta)$ is said to be {\em of general type} if $\kappa(\,(Y/\bDelta)\,) = \dim Y$. 

A Campana constellation $(X/\bDelta)$ is said to be {\em special} if there is no dominant morphism $(X/\bDelta) \to (Y/\bDelta')$ where  $(Y/\bDelta')$ is of general type.
\end{definition}

\begin{definition}
Let $f:X \to Y$ be a dominant morphism of varieties  and $(X/\bDelta)$ a Campana constellation, with $\bDelta = \sum \delta_\nu \nu$. The {\em generic fiber} of $f:  (X/\bDelta) \to Y$ is the Campana Constellation
$(X_\eta, \bDelta_\eta)$, where $X_\eta$ is the generic fiber of $f:X \to Y$, and $$\bDelta_\eta = \sum_{\nu|_{f^*\cK(Y)^\times}=0}\delta_\nu \nu,$$
namely the part of the b-divisor $\bDelta$ supported on the generic fiber.  
\end{definition}

\begin{definition} 
\begin{enumerate}
\item  Given a Campana constellation $(X/\bDelta_X)$, a dominant
  morphism $f: X \to Y$ is {\em special}, if its generic fiber is special.  
  \item[(1')] In particular, considering $X$ with trivial constellation, a dominant
  morphism $f: X \to Y$ is {\em special}, if its generic fiber is special as a variety with trivial constellation. 
  \item Given a Campana constellation $(X/\bDelta_X)$,  a proper dominant morphism $f : X \to Y$ is said to have {\em general type base} if $(Y/\bDelta_{f,\bDelta_X})$ is of general type.  
   \item[(2')] In particular, considering $X$ with trivial constellation, a proper dominant morphism $f : X \to Y$ is said to have {\em general type base} if $(Y/\bDelta_{f})$ is of general type.
\end{enumerate} 
\end{definition}

Here is the main classification theorem of Campana:

\begin{theorem}[Campana]\label{Th:Campana-core} Let $(X/\bDelta_X)$ be a Campana constellation on a projective variety $X$. There exists a dominant rational map  $c: X \das C(X)$, unique up to birational equivalence, such that 
\begin{enumerate}
\item the map $c$ has special generic fiber, and 
\item the Campana constellation base $(C(X)/\bDelta_{c,\bDelta_X})$ is of general type.
\end{enumerate} 

This map is final for (1) and initial for (2).
\end{theorem}

This is the Campana core map of $(X/\bDelta_X)$, the constellation $(C(X)/\bDelta_{c,\bDelta_X})$ being the core of $(X/\bDelta_X)$. The key case is when $X$ has the trivial constellation, and then $c: X \das (C(X)/\bDelta_{c})$ is the Campana core map of $X$ and $(C(X)/\bDelta_{c})$ the core of $X$.

\subsubsection{More examples of constellation bases} The following is a collection of examples which I find useful to keep in mind. The stated rules for the constellation bases are explained below in \ref{Sec:rule}.

\begin{enumerate}
\item Consider $f: \AA^2 \to \AA^2$ given by $s = x^2; t= y$.  We want to describe the constellation base. Clearly on $Y = \AA^2$, the divisor $\Delta_Y = 1/2 (s=0)$. But what should the multiplicity be for a divisor on some blowup of $Y$? 

The point is that $X \to Y$ is {\em toric}, and $\bDelta$ can be described using toric geometry. Indeed, the multiplicity at a divisorial valuation $\nu$ is precisely dictated by the  value of $\nu(s)$, with a simple rule: if $\nu(s)$ is even, we have $m_\nu=1$ so $\delta_\nu=0$, otherwise $m_\nu=2$ and $\delta_\nu=1/2$. 
Regular pluricanonical differentials are generated by $(ds\wedge dt)^2/s$.
\item Consider now $f: \AA^2 \to \AA^2$ given by $s = x^2; t= y^2$. The rule this time: $m_\nu=1$ and $\delta_\nu=0$ if and only if both $\nu(s)$ and $\nu(t)$ are even, otherwise  $m_\nu=2$ and $\delta_\nu=1/2$. Regular pluricanonical differentials are generated by $(ds\wedge dt)^2/st$.
\item $f: \AA^2\sqcup \AA^2 \to \AA^2$ given by $s = x_1^2; t=y_1$ and $s= x_2; t= y_2^2$
The rule this time: $m_\nu=1$ and $\delta_\nu=0$ if and only if either $\nu(s)$ or $\nu(t)$ is even, otherwise  $m_\nu=2$ and $\delta_\nu=1/2$. Regular pluricanonical differentials are generated by $ds\wedge dt$.
\item $f: X \to \AA^2$ given by the singular cover $\Spec \CC[s,t,\sqrt{st}]$. The rule: $m_\nu=1$ and $\delta_\nu=0$ if and only if either $\nu(s)+\nu(t)$ is even, otherwise  $m_\nu=2$ and $\delta_\nu=1/2$. Regular pluricanonical differentials are generated by $(ds\wedge dt)^2/st$.
\item $f: \AA^3 \to \AA^2$ given by $s = x^2y^3; t=z$. The rule: $m_\nu=1$ and $\delta_\nu=0$ if and only if either $\nu(s)=0$ or $\nu(s)\geq 2$, otherwise  $m_\nu=2$ and $\delta_\nu=1/2$. Regular pluricanonical differentials are generated by $(ds\wedge dt)^2/s$.
\end{enumerate}

\subsubsection{}\label{Sec:rule}
Where does the rule come from? When we have a toric map of affine toric varieties, we have a map of cones $f_\sigma: \sigma_X \to \sigma_Y$. Inside these cones we have lattices $N_X$ and $N_Y$ - I am considering only the part of the lattice lying in the closed cone, so it is only a monoid, not a group. The map $f_\sigma$ maps $N_X$ into a sub-monoid  $\Gamma \subset N_Y$. Each rank-1 discrete valuation $\nu$ of $Y$ has a corresponding point $n_\nu \in N_Y$, calculated by the value of $\nu$ on the monomials of $Y$: in the case of $\AA^2$ this point is simply $(\nu(s),\nu(t))$. The rule is: $m_\nu$ is the minimal positive integer such that $$m_\nu \cdot n_\nu \in \Gamma.$$

These toric examples form the basis for defining firmaments later on.

\subsubsection{Rational points and the question of integral points} 

Campana made the following bold conjecture:

\begin{conjecture}[Campana]\label{Conj:Campana-constellation}
Let $X/k$ be a variety over a number field. Then rational points are potentially dense on $X$ if and only if $X$ is special, i.e. if and only if  the core of $X$ is a point.
\end{conjecture}

It is natural to seek a good definition of integral points on a Campana constellation and translate the non-special case of the conjecture above to a conjecture on integral points on Campana constellations of general type. 

The following definition covers part of the ground. It seems natural, yet it is not satisfactory as it is quite restrictive. It is also not clear how these points behave in morphisms.  We'll be able to go a bit further with firmaments.

\begin{definition}
Let $(Y/\bDelta)$  be a Campana constellation over a number field $k$, and assume it is admissible as in Theorem \ref{Th:Campana-admissible} and discussion therein. Write as usual $\Delta_Y=\sum (1-1/m_i)\Delta_i $ for the part of $\bDelta$ with divisorial support on $Y$. Assume given a model $(\cY , \widetilde\Delta_Y)$  of $(Y,\Delta_Y)$ over $\OKS$, such that  $\cY$ is smooth and $\widetilde\Delta_Y$ a horizontal normal crossings divisor. Write $Y_0 = Y \setmin \Delta_Y$. 

Consider $y\in Y_0(K)$. We say that $y$ is a {\em soft $S$-integral point on $(Y/\bDelta)$} if for any prime $\wp$ where the Zariski closure $\bar y$ of $y$ reduces to $\widetilde\Delta$ we have $$\sum \frac{1}{m_i}\operatorname{mult}_\wp\widetilde\Delta_i\cdot \bar y \ \ \ \geq \ \ \ 1.$$ 
\end{definition}

\subsection{Bogomolov vs. Campana: some remarks about their philosophies}

Let us take a step back and reconsider what we are doing. After all, we are trying to learn something about the geometry of a variety $X$ from the data of dominant morphisms  $X \to Y$ it admits to other varieties. And somehow the effect of such a map is encoded not only in the geometry of $Y$ but in some extra structure. 

Campana's approach involves introducing a new category of objects, which I call Campana  constellations. For any dominant $f: X \to Y$, this maps leaves an indelible mark, namely a constellation, on the {\em target} $Y$, and you learn about $X$ by studying the constellations onto which it maps. 

There is an approach which is technically closely related but philosophically diametrically opposed, due to Bogomolov. Bogomolov suggests that since our object of study is $X$, we need to look for the indelible mark $f: X \to Y$ leaves {\em on $X$ itself}. Bogomolov proposes to use what has come to be called a {\em  Bogomolov sheaf}: let $d= \dim Y$ and consider the saturated image $\cF_m$ of $f^*\omega_Y^{m}$ in $\Sym^m \Omega^d_X$. These form a sheaf of algebras $\oplus_{m\geq 0} \cF_m$, and it is said to be of general type if the algebra of sections has dimension $d+1$. Bogomolov suggests that  such sheaves should have an important role in the arithmetic and geometric properties of $X$. 

Even if one prefers Bogomolov's approach, I think the achievement of Campana's Theorem \ref{Th:Campana-core} is remarkable and cannot be ignored. For example, it seems that the preprint \cite{Lu} attempted to develop a theory based entirely on Bogomolov sheaves, but the author could not resist veering towards statements such as Theorem \ref{Th:Campana-core}.

So let us take a closer look at what we have been doing with Campana's approach.

In essence, what we are trying to capture is a structure on $Y$ that measures a sort of equivalence class of dominant maps $X \to Y$. In some sense, the structure should measure to which extent the map $X \to Y$ has a section,  perhaps  locally and up to proper birational maps, or perhaps on a suitable choice of discrete valuation rings. There are some reasonable properties this should satisfy:

\begin{itemize}
\item It should be local on $Y$.
\item It should be invariant under modifications of $X$.
\item It should behave well under birational modifications of $Y$. 
\item There should be a good notion of morphism of such structures, at least on the level of dominant maps.
\end{itemize}

So far, our notion of constellation satisfies all of the above. We defined constellations in terms of divisorial valuations, which live on the function field of $Y$, and automatically behave well under birational maps. In fact I modified Campana's original definition, which relied on the divisor $\Delta_Y$, by introducing $\bDelta$ precisely for this purpose. One seems to lose in the category of computability, though not so much if one can characterise and find  admissible models. The definition was made precisely to guarantee that if $S$ is the spectrum of a complete discrete valuation ring with algebraically closed residue field, and $S \to  Y$ is {\em dominant}, then the map lifts to $S \to (Y/\bDelta)$ if, and only if, it lifts to $S \to X$.

But consider the following desirable properties, which are not yet achieved:

\begin{itemize}
\item The structure should be invariant under smooth maps on $X$.
\item In some sense it should be recovered from an open covering of $X$.
\item It should be computable.
\item There should be a notion of morphisms, good enough to work with non-dominant maps and integral points. 
\end{itemize}

It seems that Campana constellations are wonderfully suited for purposes of birational classification. Still they seem to lack some subtle information necessary to have these last properties, such as good definitions of non-dominant morphisms and integral points - at least I have not been successful in doing this directly on constellations in a satisfactory manner. For these purposes I propose the notion of firmaments. At this point I can achieve  these desired properties under extenuating circumstances, which at least enables one to state meaningful questions. It is very much possible that at the end a simpler formalism will be discovered, and the whole notion of firmaments will be redundant.

\subsection{Firmaments supporting  constellations and integral points}

The material in this section is very much incomplete as many details are missing and many questions are yet unanswered.

\subsubsection{Firmaments: valuative definition}

Let me first define the notion of firmaments in a way that seems to make things a bit more complicated than constellations, and where it is not clear that any additional desired properties are achieved. 

The underlying structure is still a datum attached to every divisorial valuation $\nu$ of $Y$. The datum is a subset $\Gamma_\nu\subset \NN$, and the sole requirement on each individual $\Gamma_\nu$ is that 
\begin{itemize}
 \item $\Gamma_\nu$ is the union of finitely many non-zero additive submonoids of $\NN$.
\end{itemize} 
and the structure is considered trivial if $\Gamma_\nu = \NN$.
 
There is an additional requirement, namely that this should come locally from a map $X \to Y$, in the way described below.

\begin{definition}
Let $Y$ be a variety, $X$ a reduced scheme, and   let $f: X \to Y$ be a  morphism, surjective on each irreducible component of $X$. For each divisorial valuation $\nu$ on $\cK(Y)$ consider $f': X_\nu' \to Y_\nu$, where $X'_\nu$ is a normal-crossings desingularization of the main part of the pullback $X \times_YY_\nu$. Write $F_\nu$ for the fiber of $X_\nu'$ over $s_\nu$.  For each point $x\in F_\nu$, assume that the components of $F_\nu$ passing through $x$ have multiplicities $m_1,\ldots,m_k$, generating a submonoid $$\Gamma_\nu^x\ \  :=\ \  \langle m_1,\ldots,m_k\rangle\ \  \subset \NN.$$
     Define 
$$\Gamma_\nu\ \  = \ \ \cup_{x\in F_\nu} \Gamma_\nu^x.$$ 
\end{definition} 

\begin{definition} A {\em firmament} $\bGamma$ on $Y$ is an assignment 
$$\nu\ \  \mapsto\ \   \Gamma_\nu\subset \NN$$
which, locally in the \'etale topology of $Y$, comes from a morphism $X \to Y$ as above.
\end{definition}

This condition on local description, which seems harmless, is actually crucial for the properties of firmaments. 

A firmament supports a unique constellation:

\begin{definition}
Let $\bGamma$ be a firmament. The {\em multiplicity} of the divisorial valuation $\nu$ is defined as $m_\nu = \min ( \Gamma_\nu \setminus \{0\})$. The {\em constellation hanging by $\bGamma$} is
$$\bDelta_\bGamma = \sum\left(1-\frac{1}{m_\nu}\right) \nu.$$
\end{definition}

\subsubsection{} Note that, according to the definition above, every firmament supports a unique constellation, though a constellation can be supported by more than one firmament. Depending on one's background, this might agree or disagree with the primitive cosmology of one's culture. Think of it this way: as we said before, the word ``constellation" refers to the entire ``heavens", visible through stronger and stronger telescopes $Y'$. The word ``firmament" refers to an overarching solid structure supporting the heavens, but solid as it may be, it is entirely imaginary and certainly not unique.  

\subsubsection{Toroidal formalism} I wish to convince the reader that this extra structure I piled on top of constellations actually makes things better. For this purpose I need to discuss a toroidal point of view. 

In fact the right foundation to use seems to be that of logarithmic structures, rather than toroidal geometry. For the longest time I stuck with toroidal geometry because the book \cite{Bourbaki} had not been written. As \cite{Ogus} is becoming available my excuses are running out, but I'll leave the translation work for the future.

\begin{definition}[\cite{KKMS}, \cite{Kato}, \cite{Abramovich-Karu}]
\begin{enumerate}
\item A toroidal embedding $U \subset X$ is the data of a variety $X$ and a dense open set $U$ with complement a Weil divisor $D= X \setminus U$,  such that locally in the \'etale, or analytic, topology, or formally, near every point,  $U\subset X$ admits an isomorphism  with (a neighborhood of a point in) $T \subset V$, with $T$ a torus and $V$ a toric variety. (It is sometimes convenient to refer to the toroidal structure using the divisor: $(X,D)$.) 
\item Let $U_X \subset X$ and $U_Y\subset Y$ be toroidal embeddings, then a dominant morphism $f:X \to Y$ is said to be toroidal if \'etale locally near every point of $X$ there is a toric chart for $X$ near $x$ and a toric chart for $Y$ near $f(x)$,  such that on these charts $f$ becomes a torus-equivariant morphism of toric varieties.
\end{enumerate}
\end{definition}

\subsubsection{The cone complex} Recall that, to a toroidal embedding $U\subset X$ we can attach  an integral polyhedral cone complex $\Sigma_X$, consisting of strictly convex cones, attached to each other along faces, and in each cone $\sigma$ a finitely generated, unit free integral saturated  monoid $N_\sigma\subset \sigma $ generating $\sigma$ as a real cone. 

Note that I am departing from usual terminology, by taking $N_\sigma$ to be the part of the lattice lying in the cone, rather than the associated group. Note also that in \cite{KKMS}, \cite{Kato} the monoid $M_\sigma$ dual to  $N_\sigma$ is used. While the use of $M_\sigma$ is natural from the point of view of logarithmic structures, all the action with firmaments happens on  its dual $N_\sigma$, so I use it instead. 

\subsubsection{Valuation rings and the cone complex} The complex $\Sigma_X$ can be pieced together using the toric charts, where the picture is well known: for a toric variety $V$, cones correspond to toric affine opens $V_\sigma$, and the lattice $N_\sigma$ is the monoid of one-parameter subgroups of the corresponding torus having a limit point in $V_\sigma$; it is dual to the lattice of effective toric Cartier divisors $M_\sigma$, which is the quotient of the lattice of regular monomials $\tilde M_\sigma$ by the unit monomials.  

For our purposes it is convenient to recall the characterization of toric cones using valuations given in \cite{KKMS}: let $R$ be a discrete valuation ring with valuation $\nu$, special point $s_R$ and generic point $\eta_R$; let $\phi:\Spec R \to X$ be a morphism such that $\phi(\eta_R) \subset U$ and  $\phi(s_R)$ lying in a stratum having chart $V=\Spec k[\tilde M_\sigma]$. One associates to $\phi$ the point $n_\phi$ in $N_\sigma$ given by the rule: 
$$n(m) = \nu (\phi^*m)\quad  \forall m\in M.$$ In case $R=R_\nu$ is a valuation ring of $Y$, I'll call this point $n_\nu$. One can indeed give a coherent picture including the case $\phi(\eta_R) \not\subset U$,  but I won't discuss this here and delay it for future treatment if one is called for. (It is however important for giving a complete picture of the category and a complete picture of the arithmetic structure.).

\subsubsection{Functoriality} Given toroidal embeddings $U_X \subset X$ and $U_Y \subset Y$ and a  morphism $f:X \to Y$ carrying $U_X$ into $U_Y$ (but not necessarily toroidal) the description above functorially associates a polyhedral morphism $f_\Sigma:\Sigma_X \to \Sigma_Y$ which is integral, that is, $f_\Sigma(N_\sigma) \subset N_\tau$ whenever $f_\Sigma(\sigma) \subset \tau$. 

\subsubsection{Toroidalizing a morphism}

While most morphisms are not toroidal, we have the following:

\begin{theorem*}[Abramovich-Karu] Let $f:X \to Y$ be a dominant morphism of varieties. Then there exist modifications $X' \to X$ and   $Y' \to Y$ and toroidal structures  $U_{X'}\subset X'$, $U_{Y'}\subset Y'$ such that the resulting rational map $f': X' \to Y'$ is a toroidal morphism:
$$\xymatrix{
U_{X'}\ar[d]\ar@{^{(}->}[r]& X'\ar[r]\ar^{f'}[d] & X\ar[d]^f \\
U_{Y'} \ar@{^{(}->}[r]& Y'\ar[r]                  & Y
}$$
 Furthermore, $f'$ can be chosen flat.
\end{theorem*}

We now define toroidal firmaments, and give an alternative  definition of firmaments in general:

\begin{definition}
A {\em toroidal firmament} on a toroidal embedding $U\subset X$ with complex $\Sigma$ is a finite collection $\bGamma= \{\Gamma_\sigma^i\subset N_\sigma\}$, where
\begin{itemize}
\item each $\Gamma_\sigma^i\subset N_\sigma$ is a finitely generate submonoid,  not-necessarily saturated.
\item each $\Gamma_\sigma^i$ generates the corresponding $\sigma$ as a cone, 
\item the collection  is closed under restrictions to faces, i.e. for each $\Gamma_\sigma^i$ and each $\tau\prec \sigma$ there is $j$ with $\Gamma_\sigma^i\cap \tau = \Gamma_\tau^j$, and
\item it is irredundant, in the sense that $\Gamma_\sigma^i \not\subset \Gamma_\sigma^j$ for different $i,j$.
\end{itemize} 

A morphism from a toroidal firmament $\bGamma_X$ on a toroidal embedding $U_X \subset X$ to   $\bGamma_Y$ on $U_Y \subset Y$ is a  morphism $f: X \to Y$ with $f(U_X)\subset U_Y$  such that for each $\sigma$ in $\Sigma_X$ and each $i$, and if $f_\Sigma(\sigma)\subset \tau$, we have  $f_\Sigma(\Gamma_\sigma^i) \subset \Gamma_\tau^j$ for some $j$. 

We say that the toroidal firmament $\bGamma_X$ is {\em induced} by $f: X \to Y$ from $\bGamma_Y$ if for each $\sigma\in \Sigma_X$ and $\tau\in \Sigma_Y$ such that $f_\Sigma (\sigma) \subset \tau$, we have  $\Gamma_\sigma^i = f_\Sigma^{-1} \Gamma_\tau^i \cap N_\sigma$.

Given a proper birational equivalence  $\phi:X_1\das X_2$, then two toroidal firmaments $\bGamma_{X_1}$ and $\bGamma_{X_2}$ are said to be {\em equivalent} if there is a toroidal embedding  $U_3 \subset X_3$, and a commutative diagram 
$$\xymatrix{
& X_3 \ar[dl]_{f_1}\ar[dr]^{f_2} \\
X_1 \ar@{-->}[rr]^\phi&&X_2,
}$$ where  $f_i$ are modifications sending $U_3$ to $U_i$, 
such that the two toroidal firmaments on $X_3$ induced by $f_i$ from  $\bGamma_{X_i}$ are identical.

A firmament on an arbitrary $X$ is the same as an equivalence class represented by  a modification $X'\to X$ with a toroidal embedding $U'\subset X'$  and a toroidal firmament $\bGamma$ on $\Sigma_{X'}$. A morphism of firmaments is a morphism of varieties which becomes a morphism of toroidal firmaments on some toroidal model.

The trivial firmament is defined by $\Gamma_\sigma = N_\sigma$ for all $\sigma$ in $\Sigma$.
\end{definition}

For the discussion below one can in fact replace $\bGamma$ by the union of the $\Gamma_\sigma^i$, but I am not convinced that makes things better.

\begin{definition} 
\begin{enumerate}
\item Let $f:X \to Y$ be a flat toroidal morphism of toroidal embeddings. The {\em base
firmament} $\bGamma_f$ associated to  $X \to Y$ is defined by the images $\Gamma_\sigma^\tau= f_\Sigma(N_\tau)$ for each cone $\tau\in \Sigma_X$ over $\sigma\in \Sigma_Y$. We make this collection irredundant by taking the sub-collection of maximal elements.
\item Let $f:X \to Y$ be a dominant  morphism of varieties. The base firmament of $f$ is represented by any $\bGamma_{f'}$, where $f': X' \to Y'$ is a flat toroidal birational model of $f$.
\item If $X$ is reducible, decomposed as $X = \cup X_i$, but $f: X_i \to Y$ is dominant for all $i$, we define the base firmament by the (maximal elements of) the union of all  the firmaments associated to $X_i \to Y$.
\end{enumerate}
\end{definition}

\subsubsection{Equivalence of definitions} Given a firmament in the new definition, given a toroidal model and given a divisorial valuation $\nu$, we have a corresponding point $n_\nu\in N_\sigma$. We define 
$$\Gamma_\nu = \{k \in \NN | kn_\nu \in \Gamma_i \text{ for some } i\}.$$
This gives a firmament in the  valuative definition.

Conversely, a firmament in the valuative definition has finitely many \'etale charts $Y_i \to Y$ where the firmament comes from $X_i \to Y_i$. One can toroidalize each $X_i \to Y_i$ simultaneously over some toroidal structure $U\subset Y$, and take the base toroidal firmament, associated $\sqcup X_i \to Y$. This gives a firmament in the ``new" sense on $Y$.

One can show that the two procedures are inverse to each other. Again I'll leave this for a later treatment.

This shows in particular that any firmament supports a unique constellation, thus allowing us access to the differential invariants of constellations. 





\subsubsection{Examples revisited}

We can now revisit our examples of base constellations in the one dimensional and higher dimensional cases, and recast them in terms of firmaments. It then becomes evident that the rules we used to calculate the constellations are simply the combinatorial data of firmaments!

\begin{enumerate}
\item $f: \AA^2 \to \AA^1$ given by $t = x^2$: $\tau = \RR_{\geq 0}; N_\tau = \NN; \bGamma=\{ 2\NN\}. $ 
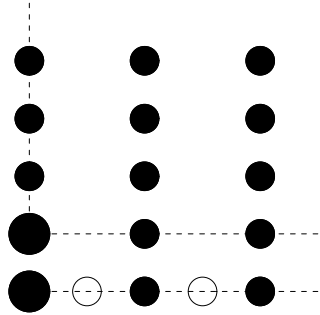
\begin{center}
\setlength{\unitlength}{0.0005in}
\begingroup\makeatletter\ifx\SetFigFont\undefined%
\gdef\SetFigFont#1#2#3#4#5{%
  \reset@font\fontsize{#1}{#2pt}%
  \fontfamily{#3}\fontseries{#4}\fontshape{#5}%
  \selectfont}%
\fi\endgroup%
{\renewcommand{\dashlinestretch}{30}
\begin{picture}(3232,3184)(0,-10)
\put(220,757){\blacken\ellipse{424}{424}}
\put(220,757){\ellipse{424}{424}}
\put(1420,757){\blacken\ellipse{300}{300}}
\put(1420,757){\ellipse{300}{300}}
\put(2620,757){\blacken\ellipse{300}{300}}
\put(2620,757){\ellipse{300}{300}}
\put(220,1357){\blacken\ellipse{300}{300}}
\put(220,1357){\ellipse{300}{300}}
\put(1420,1357){\blacken\ellipse{300}{300}}
\put(1420,1357){\ellipse{300}{300}}
\put(2620,1357){\blacken\ellipse{300}{300}}
\put(2620,1357){\ellipse{300}{300}}
\put(220,1957){\blacken\ellipse{300}{300}}
\put(220,1957){\ellipse{300}{300}}
\put(1420,1957){\blacken\ellipse{300}{300}}
\put(1420,1957){\ellipse{300}{300}}
\put(2620,1957){\blacken\ellipse{300}{300}}
\put(2620,1957){\ellipse{300}{300}}
\put(220,157){\blacken\ellipse{424}{424}}
\put(220,157){\ellipse{424}{424}}
\put(1420,157){\blacken\ellipse{300}{300}}
\put(1420,157){\ellipse{300}{300}}
\put(2620,157){\blacken\ellipse{300}{300}}
\put(2620,157){\ellipse{300}{300}}
\put(820,157){\ellipse{300}{300}}
\put(2020,157){\ellipse{300}{300}}
\put(220,2557){\blacken\ellipse{300}{300}}
\put(220,2557){\ellipse{300}{300}}
\put(1420,2557){\blacken\ellipse{300}{300}}
\put(1420,2557){\ellipse{300}{300}}
\put(2620,2557){\blacken\ellipse{300}{300}}
\put(2620,2557){\ellipse{300}{300}}
\dashline{60.000}(220,157)(3220,157)
\dashline{60.000}(220,3157)(220,757)(3220,757)
\end{picture}
}

\end{center}

\item $f: \AA^2 \to \AA^1$ given by $t = x^2y$: $\bGamma=\{ \NN\},$ the trivial structure. 
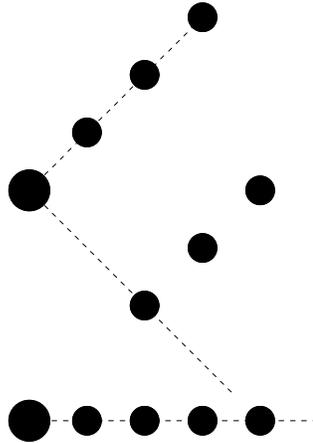
\begin{center}
\setlength{\unitlength}{0.0005in}
\begingroup\makeatletter\ifx\SetFigFont\undefined%
\gdef\SetFigFont#1#2#3#4#5{%
  \reset@font\fontsize{#1}{#2pt}%
  \fontfamily{#3}\fontseries{#4}\fontshape{#5}%
  \selectfont}%
\fi\endgroup%
{\renewcommand{\dashlinestretch}{30}
\begin{picture}(3232,4529)(0,-10)
\put(2620,2557){\blacken\ellipse{300}{300}}
\put(2620,2557){\ellipse{300}{300}}
\put(1420,3757){\blacken\ellipse{300}{300}}
\put(1420,3757){\ellipse{300}{300}}
\put(220,2557){\blacken\ellipse{424}{424}}
\put(220,2557){\ellipse{424}{424}}
\put(820,3157){\blacken\ellipse{300}{300}}
\put(820,3157){\ellipse{300}{300}}
\put(2020,4357){\blacken\ellipse{300}{300}}
\put(2020,4357){\ellipse{300}{300}}
\put(1420,1357){\blacken\ellipse{300}{300}}
\put(1420,1357){\ellipse{300}{300}}
\put(2020,1957){\blacken\ellipse{300}{300}}
\put(2020,1957){\ellipse{300}{300}}
\put(220,157){\blacken\ellipse{424}{424}}
\put(220,157){\ellipse{424}{424}}
\put(2620,157){\blacken\ellipse{300}{300}}
\put(2620,157){\ellipse{300}{300}}
\put(1420,157){\blacken\ellipse{300}{300}}
\put(1420,157){\ellipse{300}{300}}
\put(820,157){\blacken\ellipse{300}{300}}
\put(820,157){\ellipse{300}{300}}
\put(2020,157){\blacken\ellipse{300}{300}}
\put(2020,157){\ellipse{300}{300}}
\dashline{60.000}(2020,4357)(220,2557)(2320,457)
\dashline{60.000}(220,157)(3220,157)
\end{picture}
}

\end{center}
\item $f: \AA^2 \to \AA^1$ given by $t = x^2y^2$: $\bGamma=\{ 2\NN\}. $ Supported constellation: $\Delta =  D_0/2$
\item $f: \AA^2 \to \AA^1$ given by $t = x^2y^3$: $\bGamma=\{ 2\NN+ 3\NN\}. $ Supported constellation: $\Delta =  D_0/2$. Note: this is the same constellation as before, but hanging by different firmaments.
\item $f: \AA^2 \to \AA^1$ given by $t = x^3y^4$: $\bGamma=\{ 3\NN+ 4\NN\}. $ Note: even  $\Gamma\setminus \{0\}$ is not saturated in its associated group. 
\item $f: \AA^2 \to \AA^2$ given by $s = x^2; t= y$: $\bGamma=\{ 2\NN\times  \NN\}. $ 
\begin{center}
\setlength{\unitlength}{0.0005in}
\begingroup\makeatletter\ifx\SetFigFont\undefined%
\gdef\SetFigFont#1#2#3#4#5{%
  \reset@font\fontsize{#1}{#2pt}%
  \fontfamily{#3}\fontseries{#4}\fontshape{#5}%
  \selectfont}%
\fi\endgroup%
{\renewcommand{\dashlinestretch}{30}
\begin{picture}(3232,2791)(0,-10)
\put(220,219){\blacken\ellipse{424}{424}}
\put(220,219){\ellipse{424}{424}}
\put(1420,219){\blacken\ellipse{300}{300}}
\put(1420,219){\ellipse{300}{300}}
\put(2620,219){\blacken\ellipse{300}{300}}
\put(2620,219){\ellipse{300}{300}}
\put(220,819){\blacken\ellipse{300}{300}}
\put(220,819){\ellipse{300}{300}}
\put(1420,819){\blacken\ellipse{300}{300}}
\put(1420,819){\ellipse{300}{300}}
\put(2620,819){\blacken\ellipse{300}{300}}
\put(2620,819){\ellipse{300}{300}}
\put(220,1419){\blacken\ellipse{300}{300}}
\put(220,1419){\ellipse{300}{300}}
\put(1420,1419){\blacken\ellipse{300}{300}}
\put(1420,1419){\ellipse{300}{300}}
\put(2620,1419){\blacken\ellipse{300}{300}}
\put(2620,1419){\ellipse{300}{300}}
\put(220,2019){\blacken\ellipse{300}{300}}
\put(220,2019){\ellipse{300}{300}}
\put(1420,2019){\blacken\ellipse{300}{300}}
\put(1420,2019){\ellipse{300}{300}}
\put(2620,2019){\blacken\ellipse{300}{300}}
\put(2620,2019){\ellipse{300}{300}}
\put(820,219){\null\ellipse{300}{300}}
\put(820,219){\ellipse{300}{300}}
\put(2020,219){\null\ellipse{300}{300}}
\put(2020,219){\ellipse{300}{300}}
\put(820,819){\null\ellipse{300}{300}}
\put(820,819){\ellipse{300}{300}}
\put(820,1419){\null\ellipse{300}{300}}
\put(820,1419){\ellipse{300}{300}}
\put(820,2019){\null\ellipse{300}{300}}
\put(820,2019){\ellipse{300}{300}}
\put(2020,819){\null\ellipse{300}{300}}
\put(2020,819){\ellipse{300}{300}}
\put(2020,1419){\null\ellipse{300}{300}}
\put(2020,1419){\ellipse{300}{300}}
\put(2020,2019){\null\ellipse{300}{300}}
\put(2020,2019){\ellipse{300}{300}}
\put(220,2619){\blacken\ellipse{300}{300}}
\put(220,2619){\ellipse{300}{300}}
\put(820,2619){\null\ellipse{300}{300}}
\put(820,2619){\ellipse{300}{300}}
\put(1420,2619){\blacken\ellipse{300}{300}}
\put(1420,2619){\ellipse{300}{300}}
\put(2020,2619){\null\ellipse{300}{300}}
\put(2020,2619){\ellipse{300}{300}}
\put(2620,2619){\blacken\ellipse{300}{300}}
\put(2620,2619){\ellipse{300}{300}}
\dashline{60.000}(220,2619)(220,219)(3220,219)
\end{picture}
}

\end{center}
\item $f: \AA^2 \to \AA^2$ given by $s = x^2; t= y^2$: $\bGamma=\{ 2\NN\times  2\NN\}. $ 
\begin{center}
\setlength{\unitlength}{0.0005in}
\begingroup\makeatletter\ifx\SetFigFont\undefined%
\gdef\SetFigFont#1#2#3#4#5{%
  \reset@font\fontsize{#1}{#2pt}%
  \fontfamily{#3}\fontseries{#4}\fontshape{#5}%
  \selectfont}%
\fi\endgroup%
{\renewcommand{\dashlinestretch}{30}
\begin{picture}(3232,2791)(0,-10)
\put(220,219){\blacken\ellipse{424}{424}}
\put(220,219){\ellipse{424}{424}}
\put(1420,219){\blacken\ellipse{300}{300}}
\put(1420,219){\ellipse{300}{300}}
\put(2620,219){\blacken\ellipse{300}{300}}
\put(2620,219){\ellipse{300}{300}}
\put(220,1419){\blacken\ellipse{300}{300}}
\put(220,1419){\ellipse{300}{300}}
\put(1420,1419){\blacken\ellipse{300}{300}}
\put(1420,1419){\ellipse{300}{300}}
\put(2620,1419){\blacken\ellipse{300}{300}}
\put(2620,1419){\ellipse{300}{300}}
\put(820,219){\null\ellipse{300}{300}}
\put(820,219){\ellipse{300}{300}}
\put(2020,219){\null\ellipse{300}{300}}
\put(2020,219){\ellipse{300}{300}}
\put(820,819){\null\ellipse{300}{300}}
\put(820,819){\ellipse{300}{300}}
\put(820,1419){\null\ellipse{300}{300}}
\put(820,1419){\ellipse{300}{300}}
\put(820,2019){\null\ellipse{300}{300}}
\put(820,2019){\ellipse{300}{300}}
\put(2020,819){\null\ellipse{300}{300}}
\put(2020,819){\ellipse{300}{300}}
\put(2020,1419){\null\ellipse{300}{300}}
\put(2020,1419){\ellipse{300}{300}}
\put(2020,2019){\null\ellipse{300}{300}}
\put(2020,2019){\ellipse{300}{300}}
\put(220,2619){\blacken\ellipse{300}{300}}
\put(220,2619){\ellipse{300}{300}}
\put(820,2619){\null\ellipse{300}{300}}
\put(820,2619){\ellipse{300}{300}}
\put(1420,2619){\blacken\ellipse{300}{300}}
\put(1420,2619){\ellipse{300}{300}}
\put(2020,2619){\null\ellipse{300}{300}}
\put(2020,2619){\ellipse{300}{300}}
\put(2620,2619){\blacken\ellipse{300}{300}}
\put(2620,2619){\ellipse{300}{300}}
\put(220,819){\null\ellipse{300}{300}}
\put(220,819){\ellipse{300}{300}}
\put(1420,819){\null\ellipse{300}{300}}
\put(1420,819){\ellipse{300}{300}}
\put(2620,819){\null\ellipse{300}{300}}
\put(2620,819){\ellipse{300}{300}}
\put(2620,2019){\null\ellipse{300}{300}}
\put(2620,2019){\ellipse{300}{300}}
\put(1420,2019){\null\ellipse{300}{300}}
\put(1420,2019){\ellipse{300}{300}}
\put(220,2019){\null\ellipse{300}{300}}
\put(220,2019){\ellipse{300}{300}}
\dashline{60.000}(220,2619)(220,219)(3220,219)
\end{picture}
}

\end{center}
\item $f: \AA^2\sqcup \AA^2 \to \AA^2$ given by $s = x_1^2; t=y_1$ and $s= x_2; t= y_2^2$:  $\bGamma=\{ 2\NN\times  \NN, \NN\times  2\NN\}. $ Note: more than one semigroup. $\Delta_Y = 0$, but on blowup the exceptional gets $1/2$.
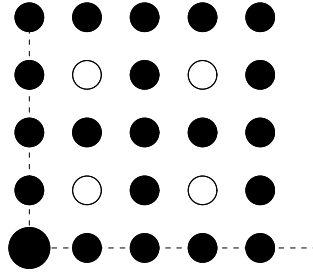
\begin{center}
\setlength{\unitlength}{0.0005in}
\begingroup\makeatletter\ifx\SetFigFont\undefined%
\gdef\SetFigFont#1#2#3#4#5{%
  \reset@font\fontsize{#1}{#2pt}%
  \fontfamily{#3}\fontseries{#4}\fontshape{#5}%
  \selectfont}%
\fi\endgroup%
{\renewcommand{\dashlinestretch}{30}
\begin{picture}(3232,2791)(0,-10)
\put(220,219){\blacken\ellipse{424}{424}}
\put(220,219){\ellipse{424}{424}}
\put(1420,219){\blacken\ellipse{300}{300}}
\put(1420,219){\ellipse{300}{300}}
\put(2620,219){\blacken\ellipse{300}{300}}
\put(2620,219){\ellipse{300}{300}}
\put(220,819){\blacken\ellipse{300}{300}}
\put(220,819){\ellipse{300}{300}}
\put(1420,819){\blacken\ellipse{300}{300}}
\put(1420,819){\ellipse{300}{300}}
\put(2620,819){\blacken\ellipse{300}{300}}
\put(2620,819){\ellipse{300}{300}}
\put(220,1419){\blacken\ellipse{300}{300}}
\put(220,1419){\ellipse{300}{300}}
\put(1420,1419){\blacken\ellipse{300}{300}}
\put(1420,1419){\ellipse{300}{300}}
\put(2620,1419){\blacken\ellipse{300}{300}}
\put(2620,1419){\ellipse{300}{300}}
\put(220,2019){\blacken\ellipse{300}{300}}
\put(220,2019){\ellipse{300}{300}}
\put(1420,2019){\blacken\ellipse{300}{300}}
\put(1420,2019){\ellipse{300}{300}}
\put(2620,2019){\blacken\ellipse{300}{300}}
\put(2620,2019){\ellipse{300}{300}}
\put(820,819){\null\ellipse{300}{300}}
\put(820,819){\ellipse{300}{300}}
\put(820,2019){\null\ellipse{300}{300}}
\put(820,2019){\ellipse{300}{300}}
\put(2020,819){\null\ellipse{300}{300}}
\put(2020,819){\ellipse{300}{300}}
\put(2020,2019){\null\ellipse{300}{300}}
\put(2020,2019){\ellipse{300}{300}}
\put(220,2619){\blacken\ellipse{300}{300}}
\put(220,2619){\ellipse{300}{300}}
\put(1420,2619){\blacken\ellipse{300}{300}}
\put(1420,2619){\ellipse{300}{300}}
\put(2620,2619){\blacken\ellipse{300}{300}}
\put(2620,2619){\ellipse{300}{300}}
\put(820,1419){\blacken\ellipse{300}{300}}
\put(820,1419){\ellipse{300}{300}}
\put(2020,1419){\blacken\ellipse{300}{300}}
\put(2020,1419){\ellipse{300}{300}}
\put(820,2619){\blacken\ellipse{300}{300}}
\put(820,2619){\ellipse{300}{300}}
\put(2020,2619){\blacken\ellipse{300}{300}}
\put(2020,2619){\ellipse{300}{300}}
\put(820,219){\blacken\ellipse{300}{300}}
\put(820,219){\ellipse{300}{300}}
\put(2020,219){\blacken\ellipse{300}{300}}
\put(2020,219){\ellipse{300}{300}}
\dashline{60.000}(220,2619)(220,219)(3220,219)
\end{picture}
}

\end{center}
\item $f: X \to \AA^2$ given by $\Spec \CC[s,t,\sqrt{st}]$: $\bGamma=\{\langle (2,0), (1,1), (0,2)\rangle \}. $ 
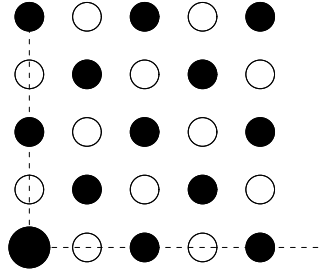
\begin{center}
\setlength{\unitlength}{0.0005in}
\begingroup\makeatletter\ifx\SetFigFont\undefined%
\gdef\SetFigFont#1#2#3#4#5{%
  \reset@font\fontsize{#1}{#2pt}%
  \fontfamily{#3}\fontseries{#4}\fontshape{#5}%
  \selectfont}%
\fi\endgroup%
{\renewcommand{\dashlinestretch}{30}
\begin{picture}(3232,2791)(0,-10)
\put(220,219){\blacken\ellipse{424}{424}}
\put(220,219){\ellipse{424}{424}}
\put(1420,219){\blacken\ellipse{300}{300}}
\put(1420,219){\ellipse{300}{300}}
\put(2620,219){\blacken\ellipse{300}{300}}
\put(2620,219){\ellipse{300}{300}}
\put(220,1419){\blacken\ellipse{300}{300}}
\put(220,1419){\ellipse{300}{300}}
\put(1420,1419){\blacken\ellipse{300}{300}}
\put(1420,1419){\ellipse{300}{300}}
\put(2620,1419){\blacken\ellipse{300}{300}}
\put(2620,1419){\ellipse{300}{300}}
\put(820,219){\null\ellipse{300}{300}}
\put(820,219){\ellipse{300}{300}}
\put(2020,219){\null\ellipse{300}{300}}
\put(2020,219){\ellipse{300}{300}}
\put(820,1419){\null\ellipse{300}{300}}
\put(820,1419){\ellipse{300}{300}}
\put(2020,1419){\null\ellipse{300}{300}}
\put(2020,1419){\ellipse{300}{300}}
\put(220,2619){\blacken\ellipse{300}{300}}
\put(220,2619){\ellipse{300}{300}}
\put(820,2619){\null\ellipse{300}{300}}
\put(820,2619){\ellipse{300}{300}}
\put(1420,2619){\blacken\ellipse{300}{300}}
\put(1420,2619){\ellipse{300}{300}}
\put(2020,2619){\null\ellipse{300}{300}}
\put(2020,2619){\ellipse{300}{300}}
\put(2620,2619){\blacken\ellipse{300}{300}}
\put(2620,2619){\ellipse{300}{300}}
\put(220,819){\null\ellipse{300}{300}}
\put(220,819){\ellipse{300}{300}}
\put(1420,819){\null\ellipse{300}{300}}
\put(1420,819){\ellipse{300}{300}}
\put(2620,819){\null\ellipse{300}{300}}
\put(2620,819){\ellipse{300}{300}}
\put(2620,2019){\null\ellipse{300}{300}}
\put(2620,2019){\ellipse{300}{300}}
\put(1420,2019){\null\ellipse{300}{300}}
\put(1420,2019){\ellipse{300}{300}}
\put(220,2019){\null\ellipse{300}{300}}
\put(220,2019){\ellipse{300}{300}}
\put(820,819){\blacken\ellipse{300}{300}}
\put(820,819){\ellipse{300}{300}}
\put(820,2019){\blacken\ellipse{300}{300}}
\put(820,2019){\ellipse{300}{300}}
\put(2020,819){\blacken\ellipse{300}{300}}
\put(2020,819){\ellipse{300}{300}}
\put(2020,2019){\blacken\ellipse{300}{300}}
\put(2020,2019){\ellipse{300}{300}}
\dashline{60.000}(220,2619)(220,219)(3220,219)
\end{picture}
}

\end{center}
\item $f: \AA^3 \to \AA^2$ given by $s = x^2y^3; t=z$: $\bGamma=\{ (2\NN + 3 \NN) \times  \NN \}. $ 
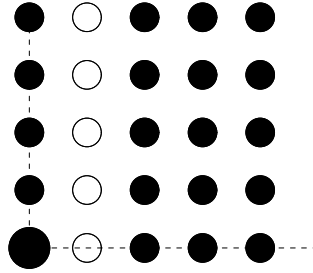
\begin{center}
\setlength{\unitlength}{0.0005in}
\begingroup\makeatletter\ifx\SetFigFont\undefined%
\gdef\SetFigFont#1#2#3#4#5{%
  \reset@font\fontsize{#1}{#2pt}%
  \fontfamily{#3}\fontseries{#4}\fontshape{#5}%
  \selectfont}%
\fi\endgroup%
{\renewcommand{\dashlinestretch}{30}
\begin{picture}(3232,2791)(0,-10)
\put(220,219){\blacken\ellipse{424}{424}}
\put(220,219){\ellipse{424}{424}}
\put(1420,219){\blacken\ellipse{300}{300}}
\put(1420,219){\ellipse{300}{300}}
\put(2620,219){\blacken\ellipse{300}{300}}
\put(2620,219){\ellipse{300}{300}}
\put(220,819){\blacken\ellipse{300}{300}}
\put(220,819){\ellipse{300}{300}}
\put(1420,819){\blacken\ellipse{300}{300}}
\put(1420,819){\ellipse{300}{300}}
\put(2620,819){\blacken\ellipse{300}{300}}
\put(2620,819){\ellipse{300}{300}}
\put(220,1419){\blacken\ellipse{300}{300}}
\put(220,1419){\ellipse{300}{300}}
\put(1420,1419){\blacken\ellipse{300}{300}}
\put(1420,1419){\ellipse{300}{300}}
\put(2620,1419){\blacken\ellipse{300}{300}}
\put(2620,1419){\ellipse{300}{300}}
\put(220,2019){\blacken\ellipse{300}{300}}
\put(220,2019){\ellipse{300}{300}}
\put(1420,2019){\blacken\ellipse{300}{300}}
\put(1420,2019){\ellipse{300}{300}}
\put(2620,2019){\blacken\ellipse{300}{300}}
\put(2620,2019){\ellipse{300}{300}}
\put(820,219){\null\ellipse{300}{300}}
\put(820,219){\ellipse{300}{300}}
\put(820,819){\null\ellipse{300}{300}}
\put(820,819){\ellipse{300}{300}}
\put(820,1419){\null\ellipse{300}{300}}
\put(820,1419){\ellipse{300}{300}}
\put(820,2019){\null\ellipse{300}{300}}
\put(820,2019){\ellipse{300}{300}}
\put(220,2619){\blacken\ellipse{300}{300}}
\put(220,2619){\ellipse{300}{300}}
\put(820,2619){\null\ellipse{300}{300}}
\put(820,2619){\ellipse{300}{300}}
\put(1420,2619){\blacken\ellipse{300}{300}}
\put(1420,2619){\ellipse{300}{300}}
\put(2620,2619){\blacken\ellipse{300}{300}}
\put(2620,2619){\ellipse{300}{300}}
\put(2020,219){\blacken\ellipse{300}{300}}
\put(2020,219){\ellipse{300}{300}}
\put(2020,1419){\blacken\ellipse{300}{300}}
\put(2020,1419){\ellipse{300}{300}}
\put(2020,2019){\blacken\ellipse{300}{300}}
\put(2020,2019){\ellipse{300}{300}}
\put(2020,2619){\blacken\ellipse{300}{300}}
\put(2020,2619){\ellipse{300}{300}}
\put(2020,819){\blacken\ellipse{300}{300}}
\put(2020,819){\ellipse{300}{300}}
\dashline{60.000}(220,2619)(220,219)(3220,219)
\end{picture}
}

\end{center}
\end{enumerate}

\subsubsection{Arithmetic} We have learned our lesson - for arithmetic we need to talk about integral points on integral models. I'll restrict to the toroidal  case, leaving the general situation to future work.

\begin{definition*} An $S$-integral model of a toroidal firmament $\bGamma$ on $Y$ consists of  an integral toroidal model  $\cY'$ of $Y'$.
\end{definition*}

\begin{definition} Consider a toroidal firmament $\bGamma$ on $Y/k$, and a rational point $y$ such that the firmament is trivial in a neighborhood of $y$. 
Let $\cY$ be a toroidal $S$-integral model.

Then $y$ is a {\em firm integral point of $\cY$} with respect to $\bGamma$ if the section $\Spec \OKS \to \cY$ is a morphism of firmaments, when $\Spec \OKS$ is endowed with the trivial firmament.

Explicitly, at each prime $\wp\in \Spec \OKS$ where $y$ reduces to a stratum with cone $\sigma$, consider the associated point $n_{y_\wp}\in N_\sigma$. Then $y$ is firmly  $S$-integral if for every $\wp$ we have
$n_{y_\wp}\in \Gamma_\sigma^i$ for some $i$.
\end{definition}

\begin{theorem} Let $f:X \to Y$ be a proper dominant morphism of varieties over $k$. 
There exists a toroidal birational model $X'\to Y'$ and an integral model $\cY'$ such that  image of a rational point on $X'$ is a firm $S$-integral point on $\cY'$ with respect to  $\bGamma_f$.
\end{theorem}

In fact, at least after throwing a few small primes into the trash-bin $S$, a point is $S$ integral on $\cY'$ with respect to  $\bGamma_f$ if and only if  locally in the \'etale topology on $\cY'$ it  lifts to a rational point on $X$. This is the motivation of the definition.

The following statements are due at least in spirit to Campana.

\begin{conjecture}
Let $(Y/\bDelta)$ be a smooth projective  Campana constellation supported by firmament $\bGamma$. Then points on $Y$ integral with respect to $\Gamma$ are potentially dense if and only if  $(Y/\bDelta)$ is special. 
\end{conjecture}

Note that this conjecture implies Conjecture \ref{Conj:Campana-constellation}:
assume this conjecture holds true.
Let $X$ be a smooth projective variety. Then rational points are potentially dense if and only if $X$ is special. 




\section{The minimal model program}
For the ``quick and easy" introduction to the minimal model program see \cite{Debarre}. For a more detailed treatment starting from surfaces see \cite{Matsuki}. For a full treatment up to 1999 see \cite{Kollar-Mori}.

The minimal model program has a beautiful beginning, a rather technical main body of work from the 80s and 90s, and quite an  exciting present. In the present account I will skip the technical main body of work.
  
\subsection{Cone of curves}
\subsubsection{Groups of divisors and curves modulo numerical equivalence} Let $X$ be a smooth complex projective variety. 

We denote by $N^1(X)$ the image of $\Pic (X) \to H^2(X, \ZZ)/\text{torsion} \subset H^2(X, \QQ)$. This is the group of Cartier divisors modulo numerical equivalence.

We denote by $N_1(X)$ the subgroup of $H_2(X,\QQ)$ generated by the fundamental classes of curves. This is the group of algebraic 1-cycles modulo numerical equivalence.

The intersection pairing restricts to $N^1(X) \times N_1(X) \to \ZZ$, which over $\QQ$ is a perfect pairing. 

\subsubsection{Cones of divisors and of curves} Denote by $Amp(X) \subset N^1(X)_\QQ$ the cone generated by classes of  ample divisors. We denote by $NEF(X)$ the closure of $Amp(X) \subset N^1(X)_\RR$, called the nef cone of $X$.

Denote by $NE(X) \subset N_1(X)_\QQ$ the cone generated by classes of  curves. We denote its closure by $\overline{NE}(X)$. The class of a curve $C$ in  $NE(X)$ is denoted $[C]$.

\begin{theorem}[Kleiman] The class $[D]$ of a Cartier divisor is in the closed cone $NEF(X)$ if and only if $[D]\cdot [C] \geq 0$ for every algebraic curve $C\subset X$. 

In other words, the cones $\overline{NE}(X)$ and $NEF(X)$ are dual to each other.
\end{theorem}

\subsection{Bend and break}

For any divisor $D$ on $X$ which is not numerically equivalent to 0, the subset $$(D \leq 0) := \{v\in NE(X) | v\cdot D \leq 0\}$$ is a half-space. The minimal model program starts with the observation that this set is especially important when $D = K_X$. In fact, in the case of surfaces, $(K_X \leq 0)\cap NE(X)$ is a subcone generated by $(-1)$-curves, which suggests that it must say something in higher dimensions. Indeed, as it turns out, it is in general a nice cone generated by so called ``extremal rays", represented by rational curves $[C]$ which can be contracted in something like a $(-1)$ contraction.

Suppose again $X$ is a smooth, projective variety with $K_X$ not nef. Our first goal is to show that there is {\em some} rational curve $C$ with $K_X \cdot C < 0 $. 

The idea is to take an arbitrary curve on $X$, and to show, using deformation theory, that it has to ``move around alot" - it has so many deformations that eventually it has to break, unless it is already the rational curve we were looking for.

\subsubsection{Breaking curves}
The key to showing that a curve breaks is the following:

\begin{lemma} 
Suppose $C$ is a projective curve of genus$>0$ with a point $p\in C$, suppose $B$ is a one dimensional affine curve, $f: C \times B \to X$ a nonconstant morphism such that $\{p\}\times B \to X$ is constant. Then, in the closure of $f(C \times B) \subset X$,
 there is a rational curve passing through $f(p)$.
  \end{lemma}
  
In genus 0 a little more will be needed:   
\begin{lemma} 
Suppose $C$ is a projective curve of genus $0$ with points $p_1,p_2\in C$, suppose $B$ is a one dimensional affine curve, $f: C \times B \to X$ a  morphism such that $\{p_i\}\times B \to X$ is constant, $i=1,2$, and the image is two-dimensional. Then the class $[f(C)]\in NE(X)$  is ``reducible": there are effective curves $C_1,C_2$ passing through $p_1, p_2$ respectively, such that $[C_1]+[C_2] = [C]$.
\end{lemma}

\subsubsection{Some deformation theory}

We need to understand deformations of a map $f: C \to X$ fixing a point or two. The key is that the tangent space of the moduli space of such maps - the deformation space -  can be computed cohomologically, and the number of equations of the deformation space is also bounded cohomologically.

\begin{lemma}
The tangent space of the deformation space of $f: C \to X$ fixing points $p_1,\ldots,p_n$ is
$$H^0\left(C, f^*T_X (-\sum p_i)\right).$$
The obstructions lie in the next cohomology group:
 $$H^1\left(C, f^*T_X (-\sum p_i)\right).$$
 The dimension of the deformation space is bounded below:
\begin{align*}
\dim Def(f: C \to X, p_1,\ldots,p_n)\ \  & \geq \ \   \chi\left(C,  f^*T_X (-\sum p_i)\right)\\
     & = \ \  - (K_X \cdot C) \ + \ (1-g(C)-n) \dim X
     \end{align*} 
  
\end{lemma} 

\subsubsection{Rational curves}
Let us consider the case where $C$ is rational. Suppose we have such a rational curve inside $X$ with $-(K_X \cdot C) \geq \dim X + 2$, and we consider deformations fixing $n=2$ of its points. Then $- (K_X \cdot C) \ + \ (1-g(C)-2) \dim X = - (K_X \cdot C) -  \dim X \geq 2$. Since $C$ is inside $X$, the only ways $f:C \to X$ can deform is either by the 1-parameter group of automorphisms, or, beyond one-parameter, go outside the image of $C$, and we get an image of dimension at least 2. So the rational curve must break, and one of the resulting components $C_1$ is a curve with $-(K_X \cdot C_1)  < -(K_X \cdot C)$. 

Suppose for a moment $-K_X$ is ample, so its intersection number with an effective curve is positive. In this case the process can only stop once we have a curve $C_\infty$ with $$0 < -(K_X \cdot C_\infty)\ \  \leq\ \  \dim X + 1.$$ Note that this is optimal - the canonical line bundle on $\PP^r$ has degree $r+1$ on any line.

\subsubsection{Higher genus} If $X$ is any projective variety with $K_X$ not nef, then there is some curve $C$ with $K_X \cdot C<0$. To be able to break $C$ we need $$ - (K_X \cdot C) - g(C) \dim X \geq 1.$$ 

There is apparently a problem: the genus term may offset the positivity of $ - (K_X \cdot C)$. One might think of  replacing $C$ by a curve covering $C$, but there is again a problem: the genus increases in coverings roughly by a factor of the degree of the cover, and this offsets the increase in $ - (K_X \cdot C)$. There is one case when this does not happen, that is in characteristic p we can take  the iterated Frobenius morphism $C^{[m]} \to C$, and the genus of $C^{[m]}$ is $g(C)$.  We can apply our bound and deduce that there is a rational curve $C'$ on $X$. If $-K_X$ is ample  we also have $0< -(K_X \cdot C')\ \  \leq\ \  \dim X + 1.$

But our variety $X$ was a complex projective variety. What do we do now? We can find a smooth model $\cX$ of $X$ over some ring $R$ finitely generated over $\ZZ$, and for each maximal ideal $\wp \subset R$ the fiber $\cX_\wp$ has a rational curve on it. 

How do we deduce that there is a rational curve on the original $X$? if $-K_X$ is ample, the same is true for $-K_\cX$, and we deduce that there is a rational curve $C_\wp$ on each $\cX_\wp$ such that $0<-(K_{\cX_\wp} \cdot C_\wp)\ \  \leq\ \  \dim X + 1.$ These are parametrized by a Hilbert scheme of finite type over $R$, and therefore this Hilbert scheme has a point over $\CC$, namely there is a rational curve $C$ on $X$ with $0<-(K_X \cdot C)\ \  \leq\ \  \dim X + 1.$

In case $-K_X$ is not ample, a more delicate argument is necessary. One fixes an ample line bundle $H$ on $\cX$, and given a curve $C$ on $X$ with $-K_X \cdot C < 0$  one shows that there is a rational curve $C'$ on each $\cX_\wp$ with
$$ H \cdot C' \leq 2 \dim X \frac{H \cdot C} {-K_X \cdot C}.$$

Then one continues with a similar Hilbert scheme argument.

\subsection{Cone theorem}

Using some additional delicate arguments one proves:

\begin{theorem}[Cone theorem] 
Let $X$ be a smooth projective variety. There is a countable collection $C_i$ of rational curves on $X$ with $$ 0 < -K_X\cdot C_i \leq \dim X + 1,$$ whose classes $[C_i]$ are discrete in the half space $N_1(X)_{K_X<0}$,  
such that 
$$\overline{NE}(X) \ \ = \overline{NE}(X)_{K_X\geq 0} + \sum_i \RR_{\geq 0} \cdot [C_i].$$
\end{theorem}

The rays $\RR_{\geq 0} \cdot [C_i]$ are called extremal rays (or, more precisely, extremal $K_X$-negative rays) of $X$.  

These extremal rays have a crucial property:

\begin{theorem}[Contraction theorem]
Let $X$ be a  smooth complex projective variety and let $R = \RR_{\geq 0} \cdot [C]$ be an extremal $K_X$-negative ray. Then there is a normal projective variety $Z$ and a surjective morphism $c_R :X \to Z$ with connected fibers, unique up to unique isomorphism, such that for an irreducible curve $D\subset X$ we have $c_R(D)$ is a point  if and only if $[D] \in R$.
\end{theorem}

This map $c_R$ is defined using a base-point-free linear system on $X$ made out of a combination of an ample sheaf $H$ and $K_X$.

\subsection{The minimal model program}

If $X$ has an extremal ray which gives a contraction to a lower dimensional variety $Z$, then the fibers of $c_R$ are rationally connected and we did learn something important about the structure of $X$: it is uniruled.

Otherwise $c_R: X \to Z$ is birational, but at least we have gotten rid of one extremal ray - one piece of obstruction for  $K_X$ to be nef.

One is tempted to apply the contraction theorem repeatedly, replacing $X$ by $Z$, until we get to a variety with $K_X$ nef. There is a problem: the variety $Z$  is often singular, and the theorems apply to smooth varieties. All we can say about $Z$ is that it has somewhat mild singularities: in general it has rational singularities; if the exceptional locus has codimension 1 - the case of a so called {\em divisorial} contraction -  the variety $Z$ has so called {\em terminal} singularities. For surfaces, terminal singularities are in fact smooth, and in fact contractions of extremal rays are just $(-1)$-contraction, and we eventually are led to a minimal model. But in higher dimensions singularities do occur.

The good news is that the theorems can be extended, in roughly the same form, to varieties with terminal singularities. (The methods are very different from what we have seen and I would rather not go into them.) So as long as we only need to deal with divisorial contractions, we can continue as in the surface case.

For non-divisorial contractions - so called small contractions - we have the following recent major result of Birkar, Cascini, Hacon and McKernan:

\begin{theorem}[Flip Theorem \cite{BCHM}]
Suppose $c_R: X \to Z$ is a small extremal contraction on a variety $X$ with terminal singularities. Then there exists another small contraction $c_R^+:X^+ \to Z$ such that $X^+$ has terminal singularities and $K_{X^+} \cdot C>0$ for any curve $C$ contracted by $c_R^+$.
\end{theorem}

The transformation $X \das X^+$ is known as a {\em flip}.

The proof of this theorem goes by way of a spectacular inductive argument, where proofs of existence of minimal models for varieties of general type, finite generation of canonical rings, and finiteness of certain minimal models are intertwined.

\begin{conjecture}[Termination Conjecture]\hfill
Any sequence of flips is finite.
\end{conjecture}

This implies the following:

\begin{conjecture}[Minimal model conjecture] 
Let $X$ be a smooth projective variety. Then either $X$ is uniruled, or there is a birational modification $X \das X'$ such that $X'$ has only terminal singularities and $K_{X'}$ is nef
\end{conjecture}

Often one combines this with the following:
\begin{conjecture}[Abundance] 
Let $X$ be a projective variety with terminal singularities and $K_X$ nef. Then for some integer $m>0$, we have $H^0(X, \cO_X(mK_X))$ is base-point-free.
\end{conjecture}

The two together are sometimes named ``the good minimal model conjecture".

The result is known for varieties of general type: it follows from the recent theorem of \cite{BCHM} on finite generation of canonical rings.

As we have seen in previous lectures, this conjecture has a number of far reaching corollaries, including Iitaka's additivity conjecture and the $(-\infty)$-conjecture. 

\section{Vojta, Campana and $abc$} In  \cite{Vojta}, Paul Vojta started a speculative investigation in diophantine geometry motivated by analogy with value distribution theory. His conjectures go in the same direction as Lang's - they are concerned with bounding the set of points on a variety rather than constructing ``many" rational points. Many of the actual proofs in the subject, such as an alternative proof of Faltings's theorem, use razor-sharp tools such as Arakelov geometry. But to describe the relevant conjectures it will suffice to discuss heights from the classical ``na\"{\i}ve" point of view. The reader is encouraged to consult Hindry--Silverman \cite{Hindry-Silverman} for a user--friendly, Arakelov--Free treatment of the theory of heights (including a proof of Faltings's theorem, following Bombieri).

A crucial feature of Vojta's conjectures is that they are not concerned just with rational points, but with algebraic points of bounded degree. To account for varying fields of definition, Vojta's conjecture always has the discriminant of the field of definition of a point $P$ accounted for. 

Vojta's conjectures are thus much farther-reaching than Lang's. You might say, much more outrageous. On the other hand, working with all extensions of a bounded degree allows for enormous flexibility in using geometric constructions in the investigation of algebraic points. So, even if one is worried about the validity of the conjectures, they serve as a wonderful testing ground for our arithmetic intuition.

\subsection{Heights and related invariants}

Consider a point in projective space $P=(x_0:\ldots:x_r) \in \PP^r$, defined over some number field $k$, with set of places $\MM_k$. Define the na\"{\i}ve height of $P$ to be

$$H(P) = \prod_{v\in \MM_k} \max( \|x_0\|_v,\ldots,\|x_r\|_v).$$

Here $\|x\|_v = |x|$ for a real $v$, $\|x\|_v = |x|^2$ for a complex $v$, and 
$\|x\|_v$ is normalized so that $\|p\| = p^{-[k_v:\QQ_p]}$ otherwise. (If the coordinates can be chosen relatively prime algebraic integers, then the product is of course a finite product over the archimedean places where everything is as easy as can be expected.)

This height is independent of the homogeneous coordinates chosen, by the product formula.

To keep things independent of a chosen field of definition, and to replace products by sums, one defines the normalized logarithmic height

$$h(P) = \frac{1}{[k:\QQ]} \log H(P).$$

Now if $X$ is a variety over $k$ with a very ample  line bundle $L$, one can consider the embedding  of $X$ in a suitable $\PP^r$ via the complete linear system of $H^0(X,L)$. We define the height $h_L(P)$ to be the height of the image point in $\PP^r$.

This definition of $h_L(P)$ is not valid for embeddings by incomplete linear systems, and is not additive in $L$. But it does satisfy these desired properties ``almost": $h_L(P) = h(P) + O(1)$ if we embed by an incomplete linear system, and $h_{L\otimes L'}(P) = h_L(P) + h_{L'}(P)$ for very ample $L, L'$. This allows us to define 

$$h_L(P) = h_A(P) - h_B(P)$$
where $A$ and  $B$ are very ample and $L\otimes B = A$. The function $h_L(P)$ is now only well defined as a function on $X(\bar k)$ up to $O(1)$.


Consider a finite set of places $S$ containing all archimedean places.

Let now $\cX$ be a scheme proper over $\OKS$, and $D$ a Cartier divisor. 

The counting function of $\cX, D$ relative to $k,S$ is a function on points of $X(\bar k)$ not lying on $D$. Suppose $P\in X(E)$, which we view again as an $S$-integral point of $\cX$. Consider a place $w$ of $E$ not lying over $S$, with residue field $\kappa(w)$. Then the restriction of $D$ to $P\simeq \Spec \cO_{E,S}$ is a fractional ideal with some multiplicity $n_w$ at $w$. We define the {\bf counting function} as follows:
$$N_{k,S}(D,P) = \frac{1}{[E:k]}\sum_{\substack{w \in \MM_E \\ w\nmid S}} n_w \log |\kappa(w)|.$$

A variant of this is the {\bf truncated counting function}

$$N^{(1)}_{k,S}(D,P) = \frac{1}{[E:k]}\sum_{\substack{w \in \MM_E \\ w\nmid S}} \min(1,n_w) \log |\kappa(w)|.$$

Counting functions and truncated counting functions depend on the choice of $S$ and a model $\cX$, but only up to $O(1)$. We'll thus suppress the subscript $S$.

One defines the {\bf relative logarithmic discriminant} of $E/k$ as follows: suppose the discriminant of a number field $k$ is denoted $D_k$. Then define 

$$d_k(E) = \frac{1}{[E:k]}\log |D_E| \ \ - \ \ \log |D_k|.$$

\subsection{Vojta's conjectures}

\begin{conjecture}
Let $X$ be a smooth proper variety over a number field $k$, $D$ a normal crossings divisor on $X$, and $A$ an ample line bundle on $X$. Let $r$ be a positive integer and $\epsilon>0$. Then there is a proper Zariski-closed subset $Z \subset X$ containing $D$ such that 
$$N_k(D,P) + d_k(k(P)) \geq h_{K_X(D)}(P)- \epsilon h_A(P) - O(1)$$
for all $P\in X(\bar k)\setminus Z$ with $[k(P):k]\leq r$.
\end{conjecture}

In the original conjecture in \cite{Vojta}, the discriminant term came with a factor $\dim X$. By the time of \cite{Vojta1} Vojta came to the conclusion that the factor was not well justified. 
A seemingly stronger version is 

\begin{conjecture}
Let $X$ be a smooth proper variety over a number field $k$, $D$ a normal crossings divisor on $X$, and $A$ an ample line bundle on $X$. Let $r$ be a positive integer and $\epsilon>0$. Then there is a proper Zariski-closed subset $Z \subset X$ containing $D$ such that 
$$N^{(1)}_k(D,P) + d_k(k(P)) \geq h_{K_X(D)}(P)- \epsilon h_A(P) - O(1).$$
\end{conjecture}
\noindent
but in \cite{Vojta1}, Vojta shows that the two conjectures are equivalent.  

\subsection{Vojta and $abc$}

The following discussion is taken from \cite{Vojta1}, section 2.

The Masser-Oesterl\'e $abc$ conjecture is the following:

\begin{conjecture}
For any $\epsilon>0$ there is $C>0$ such that for all $a,b,c\in \ZZ$, with $a+b+c=0$ and $\gcd(a,b,c) = 1$ we have
$$\max(|a|,|b|,|c|) \ \ \leq\ \ C \cdot \prod_{p|abc} p^{1+\epsilon}.$$ 
\end{conjecture}

Consider the point $P = (a:b:c) \in \PP^2$. Its height is $\log \max (|a|,|b|,|c|)$.  Of course the point lies on the line  $X$ defined by $x+y+z=0$. If we denote by $D$ the divisor of $xyz=0$, that is the intersection of $X$ with the coordinate axes, and if we set $S=\{\infty\}$, then 
$$N^{(1)}_{\QQ, S}(D,P) = \sum_{p|abc} \log p.$$

So the $abc$ conjecture says 
$$h(P) \leq (1+\epsilon) N^{(1)}_{\QQ, S}(D,P) + O(1),$$
which, writing $1-\epsilon' =  (1+\epsilon)^{-1}$, is the same as
$$(1-\epsilon') h(P) \leq  N^{(1)}_{\QQ, S}(D,P) + O(1).$$
This is applied only to rational points on $X$, so $d_\QQ(\QQ) = 0$. We have $K_X(D) = \cO_{X}(1)$, and setting $A=\cO_X(1)$ as well we get that $abc$ is equivalent to 
$$N^{(1)}_{\QQ, S}(D,P) \geq h_{K_X(D)}(P) - \epsilon'h_A(P) - O(1),$$
which is exactly what Vojta's conjecture predicts in this case.

Note that the same argument gives the $abc$ conjecture over any fixed number field.

\subsection{$abc$ and Campana}

Material in this section follows Campana's \cite{Campana2}.

Let us go back to Campana constellation curves. Recall Conjecture \ref{Conj:Cc-curve}, in particular a campana constellation curve of general type over a number field is conjectured to have a finite number of soft $S$-integral points.

Simple inequalities, along with Faltings's theorem, allow Campana to reduce to a finite number of cases, all on $\PP^1$. The multiplicities $m_i$ that occur in these ``minimal" divisors $\Delta$ on $\PP^1$ are:
$$(2,3,7), (2,4,5), (3,3,4), (2,2,2,3) \quad \text{ and }\quad (2,2,2,2,2).$$
 Now one claims that Campana conjecture in these cases follows from the $abc$ conjecture for the number field $k$. This follows from a simple application of Elkies's \cite{Elkies}. It is easiest to verify in case $k= \QQ$ when $\Delta$ is supported precisely at 3 points, with more points one needs to use Belyi maps (in the function field case one uses a proven generalization of $abc$ instead). 

We may assume $\Delta$ is supported at $0, 1$ and $\infty$.  An integral point on $(\PP^1/\Delta)$ in this case is a rational point $a/c$ such that $a,c$ are integers, satisfying the following:
\begin{itemize}
\item whenever $p | a$, in fact $p^{n_0} | a$; 
\item whenever $p | b$, in fact $p^{n_1} | b$; and
\item whenever $p | c$, in fact $p^{n_\infty} | c$,  
\end{itemize}
where $b = c-a$.

Now if $M = \max(|a|,|b|,|c|)$ then 
$$M^{1/n_0 + 1/n_1 + 1/n_\infty} \geq |a|^{1/n_0} |b|^{1/n_1} |c|^{1/n_\infty},$$
 and by assumption $a^{1/n_0} \geq \prod_{p|a}p$, and similarly for $b,c$. In other words  
$$M^{1/n_0 + 1/n_1 + 1/n_\infty} \geq \prod_{p|abc}p.$$ Since, by assumption, $1/n_0 + 1/n_1 + 1/n_\infty<1$ we can take any  $0<\epsilon<1 - 1/n_0 + 1/n_1 + 1/n_\infty$, for which the $abc$ conjecture gives  
$M^{1-\epsilon} < C\prod_{p|abc}p, $ for some $C$. So $M^{1 - 1/n_0 + 1/n_1 + 1/n_\infty-\epsilon} < C$ and $M$ is bounded, so there are only finitely many such points.

\subsection{Vojta and Campana}

I speculate: Vojta's higher dimensional conjecture implies the non-special part of Campana's conjecture \ref{Conj:Campana-constellation}, i.e. if $X$ is nonspecial its set of rational points is not dense.

The problem is precisely in understanding what happens when a point reduces to the singular locus of $D$.


\begin{thebibliography}{[20]}
\bibitem{Abramovich-quadratic} D. Abramovich, {\em  Uniformit\'e des points rationnels des courbes algŽbriques sur les extensions quadratiques et cubiques.}  C. R. Acad. Sci. Paris S\'er. I Math.  321  (1995),  no. 6, 755--758.
\bibitem{Abramovich-fibered} D. Abramovich {\em A high fibered power of a family of varieties of general type dominates a variety of general type.}  Invent. Math.  128  (1997),  no. 3, 481--494.
\bibitem{Abramovich-Karu} D. Abramovich and K. Karu, {\em Weak semistable reduction in characteristic 0. } Invent. Math.  139  (2000),  no. 2, 241--273.
\bibitem{Abramovich-Voloch} D. Abramovich and J. F. Voloch, {\em Lang's conjectures, fibered powers, and uniformity. } New York J. Math.  2  (1996), 20--34
\bibitem{Aluffi} P. Aluffi, {\em Celestial integration, stringy invariants, and Chern-Schwartz-MacPherson classes,} preprint {\tt math.AG/0506608}
\bibitem{BCHM}   C. Birkar, P. Cascini, Ch. D. Hacon and J. McKernan,
    {\em  Existence of minimal models for varieties of log general type}
preprint \url{math.AG/0610203}
\bibitem{Bourbaki} N. Bourbaki, {\em Logarithmic Structures}, Vanish and Perish  Press, Furnace Heat, Purgatory,  2015.
\bibitem{Campana-RC} F. Campana, {\em Connexit\'e rationnelle des vari\'et\'es de Fano.}\  Ann. Sci. \'Ecole Norm. Sup. (4)  25  (1992),  no. 5, 539--545. 
\bibitem{Campana} F. Campana, {\em Orbifolds, special varieties and classification theory.}  Ann. Inst. Fourier (Grenoble)  54  (2004),  no. 3, 499--630.
\bibitem{Campana2} F. Campana, {\em Fibres multiples sur les surfaces: aspects geom\'etriques, hyperboliques et arithm\'etiques. }  Manuscripta Math.  117  (2005),  no. 4, 429--461.
\bibitem{Caporaso} L. Caporaso, Lecture at MSRI, January 2006. \url{http://www.math.brown.edu/~abrmovic/GOTTINGEN/Caporaso-MSRI.pdf}
\bibitem{Caporaso-Harris-Mazur} L. Caporaso, J. Harris and B. Mazur, {\em Uniformity of rational points. } J. Amer. Math. Soc.  10  (1997),  no. 1, 1--35.
\bibitem{CT-S-SD} J.-L. Colliot-Th\'el\`ene,  A.N. Skorobogatov, and P. Swinnerton-Dyer, {\em 
Double fibres and double covers: paucity of rational points.}
Acta Arith. 79 (1997), no. 2, 113--135.
\bibitem{Corti-et-al} A. Corti et al. {\em Flips for 3-folds and 4-folds}, preprint.
\bibitem{Darmon-Granville} H. Darmon and A. Granville, {\em On the equations $z^m=F(x,y)$ and $Ax^p+By^q=Cz^r$.}  Bull. London Math. Soc.  27  (1995),  no. 6, 513--543.
\bibitem{Debarre} O. Debarre, {\em Higher-dimensional algebraic geometry. } Universitext. Springer-Verlag, New York, 2001.
\bibitem{Elkies} N. Elkies, $ABC$ implies Mordell.
Internat. Math. Res. Notices 1991, no. 7, 99--109. 
\bibitem{Faltings} G. Faltings, {\em Endlichkeitss\"atze f\"ur abelsche Variet\"aten Ÿber Zahlk\"orpern. } Invent. Math.  73  (1983),  no. 3, 349--366.
\bibitem{Fujita} T. Fujita, {\em 
On K\"ahler fiber spaces over curves.}
J. Math. Soc. Japan 30 (1978), no. 4, 779--794. 
\bibitem{Graber-Harris-Starr} T. Graber, J. Harris, and J.  Starr, {\em  Families of rationally connected varieties.}  J. Amer. Math. Soc.  16  (2003),  no. 1, 57--67 
\bibitem{Grauert} H. Grauert, {\em
Mordells Vermutung Ÿber rationale Punkte auf algebraischen Kurven und Funktionenkšrper. }
Inst. Hautes ƒtudes Sci. Publ. Math. No. 25 1965 131--149.
\bibitem{Hartshorne} R. Hartshorne, {\em  Algebraic geometry.} Graduate Texts in Mathematics, No. 52. Springer-Verlag, New York-Heidelberg, 1977.
\bibitem{Hassett} B. Hassett, {\em Correlation for surfaces of general type.}  Duke Math. J.  85  (1996),  no. 1, 95--107.
\bibitem{Hindry-Silverman} M. Hindry and J. Silverman, {\em Diophantine geometry. An introduction.} Graduate Texts in Mathematics, 201. Springer-Verlag, New York, 2000. 
\bibitem{Iitaka} S. Iitaka, {\em Algebraic geometry. An introduction to birational geometry of algebraic varieties.} Graduate Texts in Mathematics, 76. North-Holland Mathematical Library, 24. Springer-Verlag, New York-Berlin, 1982.
\bibitem{Kato} K. Kato, {\em Toric singularities. }  Amer. J. Math.  116  (1994),  no. 5, 1073--1099.
\bibitem{KKMS} G. Kempf. F. Knudsen, D. Mumford and B. Saint-Donat, {\em Toroidal embeddings. I.} Lecture Notes in Mathematics, Vol. 339. Springer-Verlag, Berlin-New York, 1973.
\bibitem{Kollar-subadditivity} J. 
Koll\'ar,
 {\em Subadditivity of the Kodaira dimension: fibers of general type. } Algebraic geometry, Sendai, 1985, 361--398,
Adv. Stud. Pure Math., 10, North-Holland, Amsterdam, 1987. 
\bibitem{Kollar} J. Koll\'ar,  {\em Rational curves on algebraic varieties. }Ergebnisse der Mathematik und ihrer Grenzgebiete. 3. Folge. A Series of Modern Surveys in Mathematics, 32. Springer-Verlag, Berlin, 1996.
\bibitem{Kollar-Miyaoka-Mori} J. Koll\'ar, Y. Miyaoka and S. Mori, {\em  Rationally connected varieties.}  J. Algebraic Geom.  1  (1992),  no. 3, 429--448.
\bibitem{Kollar-Mori} J. Koll\'ar and S. Mori, {\em Birational geometry of algebraic varieties.} With the collaboration of C. H. Clemens and A. Corti.  Cambridge Tracts in Mathematics, 134. Cambridge University Press, Cambridge, 1998. 
\bibitem{Kollar-Smith-Corti} J. Koll\'ar, K. Smith and A. Corti, {\em Rational and nearly rational varieties.} Cambridge Studies in Advanced Mathematics, 92. Cambridge University Press, Cambridge, 2004.

\bibitem{Kawamata} Y.
Kawamata, {\em 
Minimal models and the Kodaira dimension of algebraic fiber spaces.}
J. Reine Angew. Math. 363 (1985), 1--46.
\bibitem{Lang} S. Lang, {\em Fundamentals of Diophantine geometry. } Springer-Verlag, New York, 1983. 
\bibitem{Lang1} S. Lang, {\em Serge Number theory. III. Diophantine geometry.} Encyclopaedia of Mathematical Sciences, 60. Springer-Verlag, Berlin, 1991.
\bibitem{Lu}     S. S. Y. Lu,
{\em A refined Kodaira dimension and its canonical fibration} preprint \url{math.AG/0211029} 
\bibitem{Manin}
Yu, I. Manin, 
{\em 
Rational points on algebraic curves over function fields. } Izv. Akad. Nauk SSSR Ser. Mat. 27 1963 1395--1440.
Corrected in:
Yu, I. Manin, 
{\em Letter to the editors: "Rational points on algebraic curves over function fields"}
Izv. Akad. Nauk SSSR Ser. Mat. 53 (1989), no. 2, 447--448; translation in Math. USSR-Izv. 34 (1990), no. 2, 465--466

\bibitem{Matsuki} K. Matsuki, {\em Introduction to the Mori program. } Universitext. Springer-Verlag, New York, 2002.

\bibitem{Moosa-Scanlon} R. Moosa and T. Scanlon, {\em The Mordell-Lang conjecture in positive characteristic revisited.}  Model theory and applications,  273--296, Quad. Mat., 11, Aracne, Rome, 2002.
\bibitem{Ogus} A. Ogus, {\em Lectures on logarithmic algebraic geometry,} manuscript in progress 
\url{http://math.berkeley.edu/~ogus/preprints/log_book/logbook.pdf}
\bibitem{Olsson} M. Olsson, {\em  Logarithmic geometry and algebraic stacks. }  Ann. Sci. \'Ecole Norm. Sup. (4)  36  (2003),  no. 5, 747--791.
\bibitem{Pacelli} P. Pacelli, {\em Uniform boundedness for rational points. }  Duke Math. J.  88  (1997),  no. 1, 77--102.
\bibitem{Shokurov} V. V. Shokurov, {\em Prelimiting flips.}  Tr. Mat. Inst. Steklova  240  (2003),  82--219;  translation in  Proc. Steklov Inst. Math.  2003,  no. 1 (240), 75--213
\bibitem{Samuel} P. Samuel, {\em Compl\'ements \`a un article de Hans Grauert sur la conjecture de Mordell.}   Inst. Hautes \'Etudes Sci. Publ. Math. No. 29 1966 55--62.

\bibitem{Siu} Y.-T.
Siu, {\em  Invariance of plurigenera.}  Invent. Math.  134  (1998),  no. 3, 661--673.

\bibitem{Siu1} Y.-T.
Siu, {\em Extension of twisted pluricanonical sections with plurisubharmonic weight and invariance of semipositively twisted plurigenera for manifolds not necessarily of general type.}  Complex geometry (Gšttingen, 2000),  223--277, Springer, Berlin, 2002.

\bibitem{Viehweg-additivity} E. Viehweg, {\em  Die Additivit\"at der Kodaira Dimension f\"ur projektive Faserr\"aume Ÿber Variet\"aten des allgemeinen Typs. } J. Reine Angew. Math.  330  (1982), 132--142.

\bibitem{Viehweg-additivityI} E. Viehweg, {\em  Weak positivity and the additivity of the Kodaira dimension for certain fibre spaces.} Algebraic varieties and analytic varieties (Tokyo, 1981), 329--353, Adv. Stud. Pure Math., 1, North-Holland, Amsterdam, 1983. 

\bibitem{Viehweg-additivityII} E. Viehweg, {\em Weak positivity and the additivity of the Kodaira dimension. II. The local Torelli map.}  Classification of algebraic and analytic manifolds (Katata, 1982),  567--589, Progr. Math., 39, BirkhŠuser Boston, Boston, MA, 1983. 

\bibitem{Vojta} P. Vojta,  {\em  Diophantine Approximations and Value Distribution Theory,} Lecture Notes in Math. 1239, Springer, Berlin, 1987.
\bibitem{Vojta1} P. Vojta, {\em  A more general $abc$ conjecture. }  Internat. Math. Res. Notices  1998,  no. 21, 1103--1116.
\bibitem{Zariski} O. Zariski, {\em The compactness of the Riemann manifold of an abstract field of algebraic functions. }  Bull. Amer. Math. Soc.  50,  (1944). 683--691.
\end{thebibliography}
\end{document}